\documentclass[12pt,twoside]{article}
\usepackage{times}
\usepackage{amsmath,amssymb}
\usepackage{color}

\allowdisplaybreaks

\def\R {\Bbb R}
\def \N {\Bbb N}

\def\p{\partial}
\def\ve{\varepsilon}
\def\f{\frac}

\def\la{\lambda}

\def\al{\alpha}
\def\t{\tilde}
\def\q{\quad}
\def\vp{\varphi}

\def\th{\theta}
\def\g{\gamma}
\def\G{\Gamma}
\def\si{\sigma}
\def\dl{\delta}
\def\p{\partial}
\def\ve{\varepsilon}
\def\f{\frac}
\def\k{\kappa}

\def\lt{\left}
\def\rt{\right}

\def\la{\lambda}
\def\al{\alpha}
\def\wt{\widetilde}
\def\t{\tilde}
\def\r{\rho}
\def\o{\omega}

\def\vp{\varphi}
\def\th{\theta}

\def\g{\gamma}
\def\G{\Gamma}
\def\sq{\square}
\def\si{\sigma}
\def\dl{\delta}
\def\b{\beta}
\def\un{\underline}

\def\q{\qquad}

\def\ds{\displaystyle}
 \allowdisplaybreaks

\begin{document}
 \footskip=0pt
 \footnotesep=2pt
\let\oldsection\section
\renewcommand\section{\setcounter{equation}{0}\oldsection}
\renewcommand\thesection{\arabic{section}}
\renewcommand\theequation{\thesection.\arabic{equation}}
\newtheorem{claim}{\noindent Claim}[section]
\newtheorem{theorem}{\noindent Theorem}[section]
\newtheorem{lemma}{\noindent Lemma}[section]
\newtheorem{proposition}{\noindent Proposition}[section]
\newtheorem{definition}{\noindent Definition}[section]
\newtheorem{remark}{\noindent Remark}[section]
\newtheorem{corollary}{\noindent Corollary}[section]
\newtheorem{example}{\noindent Example}[section]

\title{Global smooth solutions of 3-D quasilinear wave equations with small initial data}

\author{Ding
Bingbing, \quad Liu Yingbo, \quad Yin
Huicheng\footnote{This project was supported by the NSFC (No.~11025105), and by the
Priority Academic Program Development of Jiangsu Higher Education
Institutions.}\vspace{0.5cm}\\
\small (Department of Mathematics and
IMS, Nanjing University, Nanjing 210093, China.)\\
}

\date{}
\maketitle

\centerline {\bf Abstract} \vskip 0.3 true cm

In this paper, we are concerned with the 3-D quasilinear wave equation $
\ds\sum_{i,j=0}^3g^{ij}(u, \p u)\p_{ij}^2u$ $=0$ with $(u(0,x), \p_tu(0,x))=(\ve u_0(x), \ve  u_1(x))$, where $x_0=t$,
$x=(x_1, x_2, x_3)$, $\p=(\p_0, \p_1, ..., \p_3)$, $u_0(x), u_1(x)\in C_0^\infty(\Bbb R^3)$, $\ve>0$ is small enough,
and $g^{ij}(u, \p u)=g^{ji}(u, \p u)$ are smooth in their
arguments. Without loss of generality, one can write $g^{ij}(u, \p
u)=c^{ij}+d^{ij}u+\ds\sum_{k=0}^3e^{ij}_k\p_ku+O(|u|^2+|\p u|^2)$,
where $c^{ij}, d^{ij}$ and $e^{ij}_k$ are some constants,
and $\ds\sum_{i,j=0}^3c^{ij}\p_{ij}^2=-\square\equiv -\p_t^2+\Delta$.
When $\ds\sum_{i,j,k=0}^3e^{ij}_k\o_k\o_i\o_j\not\equiv 0$ for
$\o_0=-1$ and $\o=(\o_1, \o_2, \o_3)\in\Bbb S^2$,
the authors in [7-8] have shown the blowup of the smooth solution $u$ in finite time as long as
$(u_0(x), u_1(x))\not\equiv 0$.
In the present paper, when $\ds\sum_{i,j,k=0}^3e^{ij}_k\o_k\o_i\o_j\equiv 0$,
we will prove the global existence of the smooth solution $u$. Therefore, the
complete results on the blowup or global existence of the small data solutions
have been established for the general 3-D quasilinear wave equations
$\ds\sum_{i,j=0}^3g^{ij}(u, \p u)\p_{ij}^2u=0$.

\vskip 0.3 true cm

\vskip 0.3 true cm


\vskip 0.4 true cm
\centerline{\bf $\S 1$. Introduction  and main results}
\vskip 0.3 true cm

Consider the second order quasilinear wave equation in $[0, \infty)\times \Bbb R^n$
\begin{equation}
\left\{
\begin{aligned}
&\wt\sq_g u\equiv\ds\sum_{i,j=0}^ng^{ij}(u, \p
u)\p_{ij}^2u=0,\\
&(u(0,x), \p_t u(0,x))=(\ve u_0(x), \ve u_1(x)),
\end{aligned}
\right.\tag{1.1}
\end{equation}
where $x_0=t$, $x=(x_1, ..., x_n)$, $\p=(\p_0, \p_1, ... , \p_n)$,
$\ve>0$ is a sufficiently small constant, $u_0(x), u_1(x)$ $\in
C_0^{\infty}(\Bbb R^n)$, $g^{ij}(u, \p u)=g^{ji}(u, \p u)$ are smooth functions which can be certainly expressed as
$$g^{ij}(u,\p u)=c^{ij}+d^{ij}u+\ds\sum_{k=0}^ne^{ij}_k\p_ku+O(|u|^2+|\p u|^2)\eqno{(1.2)}$$
with $c^{ij}, d^{ij}$ and $e^{ij}_k$ being some constants. Without loss of generality,
one can assume  $\ds\sum_{i,j=0}^nc^{ij}\p_{ij}^2=-\square\equiv -\p_t^2+\Delta$. By the well-known results in
[11], [15-18] and references therein, we know that (1.1) has a global smooth solution for $n\ge 4$.

If $n=3$ and $d^{ij}=0$ for all $0\le i, j\le 3$ in (1.2),
then
(1.1) has a global smooth solution if the null condition holds
(namely, $\ds\sum_{i,j,k=0}^3e^{ij}_k\o_k\o_i\o_j\equiv 0$ holds for
$\o_0=-1$ and $\o=(\o_1, \o_2, \o_3)\in\Bbb S^2$),
otherwise, the solution of (1.1) will  blow up
in finite time. See [3], [6], [10-14], [21-23] and the references therein.

If $n=3$ and $d^{ij}\not=0$ for some $(i, j)$, but  $e^{ij}_k=0$ for all
$0\le i, j, k\le 3$ in (1.2), then it follows from the
results in [4] and [19-20] that (1.1) has a global smooth
solution.

If $n=3$, $d^{ij}\not=0$ for some $(i, j)$  and  $e^{ij}_k\not=0$ for some $(i,j,k)$ in (1.2),
when $\ds\sum_{i,j,k=0}^3e^{ij}_k\o_k\o_i\o_j$ $\not\equiv 0$ for
$\o_0=-1$ and $\o=(\o_1, \o_2, \o_3)\in\Bbb S^2$,
we have established the following  blowup result
(see Remark 1.2 in [7] and Remark 1.4 in [8]): Assume $(u_0(x), u_1(x))\not\equiv 0$ and
denote by $T_{\ve}$ the lifespan of the smooth solution $u$ to (1.1), then there exists
a positive constant $\tau_0$ depending only on $(u_0(x), u_1(x))$ and the coefficients
$d^{ij}, e^{ij}_k$ in (1.2) such that
$\ds\lim_{\ve\to 0^+}\ve \ln T_{\ve}=\tau_0.$

If $n=3$, $d^{ij}\not=0$
for some $(i, j)$, and  $e^{ij}_k\not=0$ for some $(i,j,k)$ in (1.2),
when $\ds\sum_{i,j,k=0}^3e^{ij}_k\o_k\o_i\o_j$ $\equiv 0$,
a basic problem naturally arises:
does the smooth solution of
(1.1) blow up in finite time or exist globally? We will
study this problem for the equation (1.1) in the case of $n=3$
\begin{equation}
\left\{
\begin{aligned}
&\wt\sq_g u\equiv\ds\sum_{i,j=0}^3g^{ij}(u, \p
u)\p_{ij}^2u=0,\\
&(u(0,x), \p_t u(0,x))=(\ve u_0(x), \ve u_1(x)),
\end{aligned}
\right.\tag{1.3}
\end{equation}
where $u_0(x), u_1(x)\in
C_0^{\infty}(B(0, 1))$, and $B(0, 1)$ represents a unit ball centered at the origin $O$. Since the higher order terms $O(|u|^2)\p_{ij}^2u$ and $O(|\p u|^2)\p_{ij}^2u$ in (1.3)
do not play the essential roles in our study, without loss of generality, we assume that the smooth coefficients
$$g^{ij}(u, \p u)=g^{ji}(u, \p u)=c^{ij}+d^{ij}u+\ds\sum_{k=0}^3 e^{ij}_k\p_ku,\qquad i, j=0, 1, 2, 3.\eqno{(1.4)}$$
In addition,  $\ds\sum_{i,j=0}^3c^{ij}\p_{ij}^2=-\square$ is still assumed. The main conclusion in this paper is:

{\bf Theorem 1.1.} {\it Under the above assumptions on the problem (1.3), if $\ds\sum_{i,j,k=0}^3e^{ij}_k\o_k\o_i\o_j\equiv 0$
holds,
then there exists an $\ve_0>0$ such that the problem (1.3) has a global $C^\infty$ solution $u(t,x)$
for $\ve<\ve_0$.}

\vskip 0.3 true cm

{\bf Remark 1.1.} {\it The detailed decay properties of the solution $u(t,x)$ to (1.3) will be given
in (3.57)-(3.61) of Proposition 3.9 in $\S 3$ below. From Proposition 3.9, as in [4] and [19-20], we
know that the solution of (1.3) does not decay
like the solution of the 3-D free wave equation $\square v=0$ with $(v(0,x), \p_tv(0,x))=(u_0(x), u_1(x))$
since $|v|\le C(1+t)^{-1}(1+|r-t|)^{-\f12}$
holds.}

\vskip 0.2 true cm

{\bf Remark 1.2.} {\it For the 2-D nonlinear wave equation whose coefficients depend only on $\p u$
\begin{equation}
\left\{
\begin{aligned}
&\ds\sum_{i,j=0}^2g^{ij}(\p
u)\p_{ij}^2u=0,\qquad (t,x)\in [0, \infty)\times \Bbb R^2,\\
&(u(0,x), \p_t u(0,x))=(\ve u_0(x), \ve u_1(x)),
\end{aligned}
\right.\tag{1.5}
\end{equation}
where $g^{ij}(\p u)=g^{ji}(\p u)=c^{ij}+\ds\sum_{k=0}^2e^{ij}_k\p_ku
+\ds\sum_{k,l=0}^2e^{ij}_{kl}\p_ku\p_lu+O(|\p u|^3)$, $\ve>0$ is sufficiently small,
and  $\ds\sum_{i,j=0}^2c^{ij}\p_{ij}^2u=-\square$.
As $\ds\sum_{i,j,k=0}^2e^{ij}_k\o_k\o_i\o_j$ $\not\equiv0$
or $\ds\sum_{i,j,k,l=0}^2e^{ij}_{kl}\o_k\o_l\o_i\o_j\not\equiv0$ for
$\o_0=-1$ and $\o=(\o_1, \o_2)\in\Bbb S^1$,
it is well-known that the smooth solution to (1.5) must blow up in finite time
as long as $(u_0(x), u_1(x))\not\equiv0$ (see [1], [3], [9-11], [18] and so on).
As $\ds\sum_{i,j,k=0}^2e^{ij}_k\o_k\o_i\o_j\equiv0$
and $\ds\sum_{i,j,k,l=0}^2e^{ij}_{kl}\o_k\o_l\o_i\o_j\equiv0$,
(1.5) has a global smooth solution (see [2]).  In summary, the
complete results on the blowup or global existence have been established for the small data smooth solution  problem (1.5).}

\vskip 0.2 true cm

{\bf Remark 1.3.} {\it Consider  more general 2-D nonlinear wave equations whose coefficients
depend on the solution $u$ as well as its first order derivatives $\p u$
\begin{equation}
\left\{
\begin{aligned}
&\ds\sum_{i,j=0}^2g^{ij}(u, \p
u)\p_{ij}^2u=0,\qquad (t,x)\in [0, \infty)\times \Bbb R^2,\\
&(u(0,x), \p_t u(0,x))=(\ve u_0(x), \ve u_1(x)),
\end{aligned}
\right.\tag{1.6}
\end{equation}
where $g^{ij}(u,\p u)=g^{ji}(u,\p u)=c^{ij}+d^{ij}u+\ds\sum_{k=0}^2e^{ij}_k\p_ku+\t d^{ij}u^2+\ds\sum_{k=0}^2\t e^{ij}_ku\p_ku
+\ds\sum_{k=0}^2e^{ij}_{kl}\p_ku\p_lu+O(|u|^3+|\p u|^3)$, $d^{ij}\not=0$ for some $(i,j)$,
$\ve>0$ is sufficiently small, and  $\ds\sum_{i,j=0}^2c^{ij}\p_{ij}^2u=-\square$. When $\ds\sum_{i,j,k=0}^2e^{ij}_k\o_k\o_i\o_j\not\equiv0$
and $(u_0(x), u_1(x))\not\equiv 0$,
we have shown the blowup of the smooth solution $u$ to (1.6) in finite time and further given a precise description
of the blowup mechanism in [7]. With respect to the other left
cases for the coefficients $d^{ij}, e^{ij}_k, \t d^{ij}, \t e^{ij}_k$  and $e^{ij}_{kl}$ of $g^{ij}(u, \p u)$,
by our knowledge   so far there are no complete results
on the blowup or global existence of the solution $u$ to (1.6).}

\vskip 0.2 true cm

We now give the comments on the proof of Theorem 1.1. To prove Theorem 1.1, we will use the usual energy method
and the continuity argument. For this purpose, we define the standard energy
$E_N(t)=\ds\sum_{0\le I\le N}\int |\p Z^Iu|^2dx$ of (1.3)
for some suitably large integer $N$, where $Z$ denotes one of the Klainerman's fields
$\{\p_t, \p_\al, S=t\p_t+\ds\sum_{i=1}^3x_i\p_i, \G_{0\al}=t\p_\al+x_\al\p_t,
\G_{\al\b}=x_\b\p_\al-x_\al\p_\b, \al,\b=1,2,3\}$.
Motivated by [4] and [19-20], where the global existence of small data solutions
to the 3-D quasilinear wave equations
$\ds\sum_{i,j=0}^3g^{ij}(u)\p_{ij}^2u=0$ is established, we think that the solution of (1.3) does not decay
like the solution of the free linear wave equation even if the null condition (i.e, $\ds\sum_{i,j,k=0}^3e^{ij}_k\o_k\o_i\o_j\equiv 0$
for $\o_0=-1$ and $\o=(\o_1, \o_2, \o_3)\in\Bbb S^2$) holds. To prove the
global existence of the solution to (1.3), as in [20], at first we assume that $E_N(t)\le C\ve^2(1+t)^{\dl}$ holds for $t\in [0, T]$
and some fixed constant $\dl$ with $0<\dl<\ds\f12$, then we manage to derive a strong energy estimate
$E_N(t)\le C\ve^2(1+t)^{C\ve}$ for sufficiently small  $\ve>0$. In this process, we will adopt the
energy method with some special weight depending on the solution of an approximate
eikonal equation of (1.3) (see (2.34) in $\S 2$ below) as well as the property of the null
condition $\ds\sum_{i,j,k=0}^3e^{ij}_k\o_k\o_i\o_j\equiv 0$.
Here we point out that
although some main procedures in this paper are analogous to those in [20] for considering the 3-D
wave equations $\ds\sum_{i,j=0}^3g^{ij}(u)\p_{ij}^2u=0$, our analysis
is more involved since the coefficients of (1.3) depend on $u$ and $\p u$ simultaneously, and meanwhile
the null condition property of (1.3)
should be specially paid attention. Finally, based on the apriori energy estimates mentioned above, we can complete the proof
of Theorem 1.1 by the local existence of the solution to (1.3) and the continuity argument.

The rest of the paper is organized as follows. In $\S 2$, we will
give some preliminary
knowledge on the null frame in Lorentzian metric $g^{ij}(u,\p u)$ and some properties of
null conditions, where the null frames in Lorentzian metric $g^{ij}(u)$ are introduced and applied in [5] and [20].
In addition, some useful estimates and calculations are listed
or proved. In $\S 3$, we will derive the sharp decay estimate of the solution $u$ to (1.3) under the assumption
of the weak decay estimate $|Z^I u|\le C\ve (1+t)^{-\nu}$ for $1/2<\nu<1 $ and $0\le I\le N-3 $.
In $\S 4$, under some suitable assumptions on the solution $u$ to (1.3), the global weighted energy estimate on  $u$
is established by choosing an appropriate weight
so that the null condition can be utilized sufficiently. In $\S 5$, by establishing a Poincar\'{e}-type lemma similar to Lemma 8.1 of [20],
we can obtain the higher order energy estimate on the solution $u$ to (1.3). Finally,
in $\S 6$, we will derive the precise energy estimate $E_N(t)\le C\ve^2(1+t)^{C\ve}$.
Meanwhile, all the assumptions on $u$ in $\S 3$-$\S 5$ are closed. Therefore,
the proof of Theorem 1.1 is completed by continuity argument.

\vskip 0.2 true cm

In the whole paper, we will use the following notations:

Let $Z$ denote one of the Klainerman's fields
$$
\p_t, \p_\al, S=t\p_t+\ds\sum_{i=1}^3x_i\p_i, \G_{0\al}=t\p_\al+x_\al\p_t,
\G_{\al\b}=x_\b\p_\al-x_\al\p_\b, \al,\b=1, 2, 3.
$$

Let $\p$ stand for $\p_t$ or $\p_\al$ ($\al=1,2,3$).
The norm $\|f\|_{L^2}$ means $\|f(t,\cdot)\|_{L^2(\Bbb R^3)}$.
$C$ denotes by a generic positive constant, and
$|Z^kv|\equiv \ds\sum_{|\nu|=k}|Z^{\nu}v|$ for $k\in\Bbb N\cup\{0\}$(or bold {\bf k}, $I, J, K\in\Bbb N\cup\{0\}$)
and the multiple indices $\nu$'s. Specially, we denote $|Zu|$ by $|Z^1 u|$.

\vskip 0.4 true cm
\centerline{\bf $\S 2$. Preliminaries}

\vskip 0.4 true cm

As in [5] and [20], we introduce the following nullframe $\{L,\ds\underline {L}\;, S_1,S_2\}$ for the Minkowski metric
$ds^2=-dx_0^2+\ds\sum_{\al=1}^3dx_\al^2$:
\begin{equation}
\left\{
\begin{aligned}
&L=L^i\p_i,\;\qquad \;\underline L=\underline L^i\p_i,\\
&S_1=\ds\f{1}{\sqrt {1-(\o^3)^2}}(\o^3\o^\al-\dl^{\al3})\p_\al,\qquad \;S_2=\ds\f{1}{\sqrt{1-(\o^3)^2}}(-\o^2\p_1+\o^1\p_2),
\end{aligned}
\right.\tag{2.1}
\end{equation}
here and in what follows the repeated upper and lower indices stand for the summations over $i=0,1,2,3$ (or $j,k,m=0,1,2,3$)
and $\al=1,2,3$ (or $\b=1,2,3$) respectively, and
$$L^0={\un L}^0=1,\; L^{\al}=\o ^\al,\; \un L ^\al=-\o^\al,\; \o^\al=\f{x_\al}{|x|},\; \al=1,2,3.$$
In addition, one can raise and lower the indices with respect to the Minkowski metric $ds^2=-dx_0^2+\ds\sum_{\al=1}^3dx_\al^2
\equiv\ds\sum_{i,j=0}^3c_{ij}dx_idx_j$
as follows:  $X_i=c_{ij}X^j, \;X^j=c^{ij}X_j $,
where the matrix $(c^{ij})_{i,j=0}^3$
stands for the inverse of the
matrix $(c_{ij})_{i,j=0}^3$ (more concretely,
$c^{00}=-1, c^{\al\al}=1$, and $c^{ij}=0$  for $i\not=j$, which coincides with
the assumption of $c^{ij}\p_{ij}^2=-\square$ for the equation (1.3)).

Since $\{ L,\ds\underline {L}\;, S_1,S_2\}$ is a null frame, we can make a decomposition for the Lorentian metric $g^{ij}$
in the equation (1.3) (here one notes that $g^{ij}(u,\p u)$ in (1.3) is somewhat different from $g^{ij}(u)$ in [20]):
$$g^{ij}=g^{UV}U^iV^j,\eqno{(2.2)}$$
where $U,V \in\{L,\underline {L}\;,S_1, S_2\}$. As in (2.7) of [20], we define
$$g_{UV}=g^{ij}U_iV_j,\eqno{(2.3)}$$
which denotes the lowering of the indices and not the inverse of $g^{UV}$.
By (2.6) of [20], $g^{UV}$  and  $g_{UV}$  have the following relations:
$$ g^{L\underline {L}}=\ds\f{1}{4}g_{\underline {L}L},\;g^{LL}=\ds\f{1}{4}g_{\underline L\;\underline L},\;g^{\underline L\;\underline L}=\ds\f{1}{4}g_{LL},\;
g^{LA}=-\ds\f{1}{2}g_{\underline {L} A},\;g^{\underline {L}A}=-\ds\f{1}{2}g_{LA},\;g^{AB}=g_{AB},\eqno{(2.4)}$$
here and below $A, B$ denote any of the vectors $S_1, S_2$.

Through the whole paper, we denote by
$$\bar\p=\{L ,S_1, S_2\},\quad \text{which are tangent to the outgoing light cone $\{x_0^2=\ds\sum_{\al=1}^3x_\al^2\}$},$$
and
$$ \p_q=\ds\f{1}{2}(\p_r-\p_t),\;\;\p_p=\ds\f{1}{2}(\p_r+\p_t).$$
Based on the preparations above, next we cite or establish some inequalities  which will be used frequently.

{\bf Lemma 2.1. (See [20] Lemma 3.1)} {\it For any smooth function $u$
$$(1+t+r)|\bar\p u|\leq C |Z u|,\eqno{(2.5)}$$
$$(1+|t-r|)^k|\p^k u|\leq C\ds\sum_{0\le I\leq k}|Z^I u|, \eqno{(2.6)}$$
$$ (1+t+r)|\p u|\leq Cr|\p_q u|+C|Z u|,\eqno{(2.7)}$$
$$|\bar\p^2 u|+r^{-1}|\bar\p u|\leq \ds\f{C} {r}\ds\sum_{0\le I\leq 2}\ds\f{|Z^I u|}{1+t+r},\quad\text {where $|\bar\p^2 u|^2=\ds\sum_{S,T\in \Bbb T}|ST u|^2$,\quad $\Bbb T=\{L, S_1, S_2\}$}.\eqno{(2.8)}$$}

Next, as in Lemma 3.2 of [20], we look for some ``good derivatives'' so that
$\wt\sq_g u$ in (1.3) can be approximated well. To this end, we set
$$ D^{ij}=g^{ij}-c^{ij}-E^{ij},\eqno{(2.9)}$$
where
$$ E^{ij}=\ds\sum^3_{k=0}e^{ij}_k\p_k u. \eqno{(2.10)}$$

Let $L^i_1=-\ds\f{1}{2}g_{L\un L}L^i-\ds\f{1}{4}g_{LL}\un {L}^i+g_{LA} A^i$ (this definition is the same as that in (2.9) of [20, Lemma 2.1]),
where $g_{LA} A^i=g_{LS_1} S_1^i+g_{LS_2} S_2^i$.
Together with (2.3) and (2.9), this yields
$$ L^i_1=L^i-\ds\f{1}{2}D_{L\un L}L^i-\ds\f{1}{4}D_{LL}\un{L}^i+D_{LA} A^i-\ds\f{1}{2}E_{L\un L}L^i-\ds\f{1}{4}E_{LL}\un{L}^i+E_{LA} A^i
\eqno{(2.11)}$$
In addition, we introduce another somewhat different $L^i_2$ from that in (3.8) of [20] due to the appearance of $E^{ij}$ in the coefficients $g^{ij}(u,\p u)$:
$$ L^i_2=L^i-\ds\f{1}{4}D_{LL}\un {L}^i-\ds\f{1}{4}E_{LL}\un {L}^i. \eqno{(2.12)}$$
Similar to Lemma 3.2 of [20], we have

{\bf Lemma 2.2.} {\it Assume  $|D|\leq \ds\f{1}{4}$ and $|E|\leq \ds\f{1}{8}$, then
$$ |(2L^i_1\p_i+\ds\f{\ell}{r})(r\p_q u)-r\wt\sq _g u|\leq \ds\f{C}{1+t+r}\ds\sum_{0\le I\leq 2}|Z^I u|,\eqno{(2.13)}$$
where $\ell=\bar{tr}D+D_{L\un L}-\ds\f{1}{2}D_{LL}+\bar{tr}E+E_{L\un L}-\ds\f{1}{2}E_{LL}$ with $\bar{tr}D=\dl^{AB}D_{AB}$,
$\bar{tr}E=\dl^{AB}E_{AB}$ and $\dl^{AB}$ the Kronecher delta function. If assume also that
$$|D_{LL}|+|D_{LA}+|D_{AA}|+|D_{L\un L}| \leq \ds\f{1+|t-r|}{1+t+r}\lt|\ds\f{1+t+r}{1+|t-r|}\rt|^a \eqno{(2.14)}$$
and
$$ |E_{LL}|+|E_{LA}|+|E_{AA}|+|E_{L\un L}| \leq \ds\f{1}{1+t+r}\lt|\ds\f{1+t+r}{1+|t-r|}\rt|^a ,\;\;a\ge 0,\eqno{(2.15)}$$
then
$$|2L^i_2\p_i(r\p_q u)-r\wt\sq_g u|\leq \ds\f{C}{1+t+r}\lt|\ds\f{1+t+r}{1+|t-r|}\rt|^a \ds\sum_{0\le I\leq 2}|Z^I u|.\eqno{(2.16)}$$}
{\bf Proof.}  At first, we point out that although the conclusions in Lemma 2.2 are rather analogous to the ones in Lemma 3.2 of [20],
we still give the detailed proof since the coefficients $g^{ij}(u,\p u)$ and the operator $\ell$ in (2.13)
are somewhat different from the corresponding ones in [20].

Note that
$$\bar{tr}g=\dl^{AB}g_{AB}=2+\bar{tr}D+\bar{tr}E\quad \text{and} \quad 2L^i_1\p_i r=2-D_{L\un L}+\ds\f{1}{2}D_{LL}-E_{L\un L}+\ds\f{1}{2}E_{LL},$$
which derives $\ell=\bar{tr}g-2L^i_1\p_i r$. Then it follows from a direct computation that
\begin{align*}
 &(2L^i_1\p_i+\ds\f{\ell}{r})(r\p_q u)-r\wt\sq _g u\\
 =&2L^i_1\p_i(r\p_q u)-2(L^i_1\p_i r)\p_q u+(\bar{tr}g)\p_q u-r\wt\sq_g u\\
 =&r\bigl(2L^i_1\p_i\p_q u+\ds\f{\bar tr g}{r}\p_q u-\wt\sq_g u\bigr). \tag{2.17}
\end{align*}
By $|2L^i_1\p_i\p_q u+\ds\f{\bar trg}{r}\p_q u-\wt\sq_g u|\leq C\ds\sum_{1\leq k\leq 2}r^{k-2}|\bar\p^k u|$ (
see (2.16) of [20, Lemma 2.2]) together with (2.8), thus we have
$$
|(2L^i_1\p_i+\ds\f{\ell}{r})(r\p_q u)-r\wt\sq _g u|
\le C\ds\sum_{0\le I\le 2}\ds\f{|Z^I u|}{1+t+r},\eqno{(2.18)}
$$
which means that (2.13) holds.

In addition, it follows from $L^i\p_i r=1$ and $A^i\p_i r=0$ that
\begin{align*}
&2L^i_2\p_i(r\p_q u)-r\wt\sq_g u\\
&=2L^i_1\p_i(r\p_q u)+2(\ds\f{1}{2}D_{L\un L}L^i-D_{LA} A^i+\ds\f{1}{2}E_{L\un L}L^i-E_{LA} A^i)\p_i(r\p_q u)-r\wt\sq_g u\\
&\equiv I+II+III,\tag{2.19}
\end{align*}
where
\begin{align*}
&I=(2L^i_1\p_i+\ds\f{\ell}{r})(r\p_q u)-r\wt\sq_g u,\\
&II=-(\bar tr D-\ds\f{1}{2}D_{LL}+\bar tr E-\ds\f{1}{2}E_{LL})\p_q u\\
&III=r(D_{L\un L}+E_{L\un L})L^i\p_i\p_q u-2r(D_{LA}+E_{LA})A^i\p_i\p_q u.
\end{align*}
Since the first term $I$ has been estimated in (2.18), we only need to treat $II$ and $III$
in (2.19). By (2.14)-(2.15) and (2.6), one easily knows
$$
|II|\le C\ds\f{1+|t-r|}{1+t+r}\lt|\ds\f{1+t+r}{1+|t-r|}\rt|^a|\p u|
\le\ds\f{C}{1+t+r}\lt|\ds\f{1+t+r}{1+|t-r|}\rt|^a|Z u|.\eqno{(2.20)}
$$
On the other hand, it also follows from $L^i\p_i \o^j=0$ and (2.5)-(2.6) that
\begin{align*}
&r|L^i\p _i(\un L^j\p_j u)|=r|L^i\un L^j\p_i\p_j u|\le Cr|\bar \p \p u|\le C|Z \p u|\le C\ds\sum_{0\le I\le 1}|\p Z^I u|\\
&\le \ds\f{C}{1+|t-r|}\ds\sum_{0\le I\le 2}|Z^I u|,
\end{align*}
which derives
$$r|L^i\p_i\p_q u|=\f{r}{2}|L^i\p _i(\un L^j\p_j u)|\le\ds\f{C}{1+|t-r|}\ds\sum_{0\le I\le 2}|Z^I u|.\eqno{(2.21)}$$
Similarly, one has
$$r|A^i\p_i\p_q u|\le\ds\f{C}{1+|t-r|}\ds\sum_{0\le I\le 2}|Z^I u|.\eqno{(2.22)}$$
Thus, collecting (2.21)-(2.22) together with (2.14)-(2.15) yields
\begin{align*}
|III|&\le \ds\f{C}{1+|t-r|}(|D_{L\un L}|+|E_{L\un L}|+|D_{LA}|+|E_{LA}|)\ds\sum_{0\le I\le 2}|Z^I u|\\
&\le \ds\f{C}{1+t+r}\lt|\ds\f{1+t+r}{1+|t-r|}\rt|^a \ds\sum_{0\le I\leq 2}|Z^I u|.\tag{2.23}
\end{align*}
Finally, substituting (2.18), (2.20) and (2.23) into (2.19) yields (2.16). \qquad $\square$

\vskip 0.2 true cm
As in [20] and [4], in order to obtain the sharp decay of the solution $u$ to (1.3),
we will adopt the idea of  integration along the integral curves for the eikonal equation of (1.3)
(see (2.34) below)
so that the usual phase $r-t$ can be replaced and further the decay estimates of $u$
on a curved background can be treated conveniently.
\vskip 0.3 true cm
Let $X_{\la}(s)=(X_{\la}^0(s), X_{\la}^1(s), X_{\la}^2(s), X_{\la}^3(s))\in \Bbb R^{1+3}$ ($\la=1,2)$  be the backward integral curve of the vector
fields $ L_\la=L_{\la}^i\p_i$. Namely, the components of $X_{\la}(s)$ satisfy
\begin{equation}
\left\{
\begin{aligned}
&\ds\f{d}{ds}X^i_{\la}(s)=L^i_\la(X_\la(s)), \;\;s\leq 0,\\
&X_{\la}(0)=(t,x),
\end{aligned}
\right.\tag{2.24}
\end{equation}
where  $(t,x)\in H\equiv\{(t,x)\in \Bbb R^+\times\Bbb R^3: \ds\f{t}{2}<|x|<\f{3t}{2}\}$  which is near the light cone.
Let $s_\la<0$ be the largest number such that $ X_\la (s_\la)\in \p H$ and $X_\la (s)\in H$ for $s>s_\la$, where $\p H$
represents the boundary of $H$.
Define $ \tau_\la=\tau_\la(t,x)=X^0_\la (s_\la)$. When we assume $|D|\le 1/4$ and $|E|\le 1/8$,  then  the
corresponding integral curves will intersect $\p H$ obviously. Next we establish some results similar to
the ones in Lemma 4.1 of [20] so that the decay estimates of the first order derivatives $\p u$  can be
controlled.

{\bf Lemma 2.3.} {\it Assume $ |D|\leq 1/16 ,\;|E|\leq 1/32$ and
$$\int^T_0 (\|D(t,.)\|_{L^\infty(H_t)}+\|E(t,.)\|_{L^\infty(H_t)})\ds\f{dt}{1+t}\leq 1,\eqno{(2.25)}$$
where $H_t\equiv\{x\in\Bbb R^3: t/2<|x|<3t/2\}$, then
\begin{align*}
&(1+t+r)|\p u(t,x)|\leq C \sup\limits_{\tau_1\leq\tau\le t}\ds\sum_{0\le I\leq 1}\|Z^I u(\tau,.)\|_{L^\infty}\\
&\qquad +C\int^t_{\tau_1}\biggl((1+\tau)\|\wt\sq_g u(\tau,.)\|_{L^\infty(H_\tau)}+\ds\sum_{0\le I\leq 2}(1+\tau)^{-1}\|Z^I u(\tau,.)\|_{L^\infty(H_\tau)}\biggr)d\tau.\tag{2.26}
\end{align*}

If we assume $ |D|\leq 1/16 ,\;|E|\leq 1/32$, and for some constant $a\ge 0$
$$|D_{LL}|+|D_{LA}|+|D_{AA}|+|D_{L\un L}|\le \ds\f{1}{4}\ds\f{1+|t-r|}{1+t+r}\lt|\ds\f{1+t+r}{1+|t-r|}\rt|^a  \;\;in\quad H, $$
$$|E_{LL}|+|E_{LA}|+|E_{AA}|+|E_{L\un L}| \leq \ds\f{1}{4}\ds\f{1}{1+t+r}\lt|\ds\f{1+t+r}{1+|t-r|}\rt|^a,\eqno{(2.27)}$$
then
\begin{align*}
&(1+t+r)|\p u(t,x)|\leq C \sup\limits_{\tau_2\leq\tau\le t}\ds\sum_{0\le I\leq 1}\|Z^I u(\tau,.)\|_{L^\infty}\\
&\quad +C\int^t_{\tau_2}\biggl((1+\tau)\|\wt\sq_g u(\tau,.)\|_{L^\infty(H_\tau)}+\ds\sum_{0\le I\leq 2}(1+\tau)^{-1+a}\|(1+|q(\tau,.)|)^{-a}Z^I u(\tau,.)\|_{L^\infty(H_\tau)}\biggr) d\tau,\tag{2.28}
\end{align*}
where $q(t,x)=r-t$.}

\vskip 0.2 true cm
{\bf Remark 2.1.} {\it By comparison with the assumptions in Lemma 4.1 of [20], we have posed an extra
decay assumption on the term $|E_{LL}|+|E_{LA}|+|E_{AA}|+|E_{L\un L}|$ in (2.27). Notice that there is no such an
assumption (2.27) in [20] since there coefficients $g^{ij}(u)$ depend only on the solution itself.}
\vskip 0.2 true cm

{\bf Proof.}  By (2.7) it suffices to prove that $\phi=r\p_q u$ can be controlled by the right hand side
of (2.26) or (2.28). To this end, we will divide the proof process into the following two cases
of $(t,x)\not\in H$ and $(t,x)\in H$ separately.

{\bf Case A.\quad $(t,x)\not\in H$}

In this case, one has $|t-r|\ge t/2$. This, together with $r\le 1+t$ and  ${\un L}=\ds\f{S-\o^\al\G_{0\al}}{r-t}$, yields
$$|\phi(t,x)|=\ds\f12|r{\un L}u|\le C(1+2|t-r|)|{\un L} u|\le C\ds\sum_{0\le I\leq 1}\|Z^I u(\tau,.)\|_{L^\infty}.\eqno{(2.29)}$$

{\bf Case B.\quad $(t,x)\in H$}

By the characteristics method and (2.24), we have
$$\ds\f{d}{ds}\left(\phi(X_1(s))e^{G_1(s)} \right)=\ds\f{1}{2}e^{G_1(s)}\lt((2L^i_1\p_i+\ds\f{\ell}{r})\phi\rt)(X_1(s)),\eqno{(2.30)}$$
where $G_1(s)=-\ds\f{1}{2}\int^0_s \ds\f{\ell(X_1(\si))}{r(X_1(\si))}d\si$, and the quantity $\ell$
has been defined in (2.13) of Lemma 2.2.

If $\tau_1=\tau_1(t,x)\ge 1$, then $\tau\ge \tau_1 \ge 1$ for $s\ge s_1$ and further
 $r\ge \ds\f{\tau}{2}\ge \ds\f{1+\tau}{4}$ for $(\tau, x) \in X_1(s)\in H$. Thus, one has  under the assumption (2.25)
\begin{align*}
|G_1(s)|&\le \ds\f{1}{2}\int^0_s|\ds\f{\ell}{r}(X_1(\si))|d\si\le \ds\f{1}{2}\int^0_{s_1}|\ds\f{\ell}{r}(X_1(\si))|d\si\\
&\le C \int^t_{\tau_1} \ds\f{\|D(\tau,.)\|_{L^\infty(H_\tau)}+\|E(\tau,.)\|_{L^\infty(H_\tau)}}{1+\tau}d\tau\le C.\tag{2.31}
\end{align*}

In addition, by $t=X^0_1(s)$ and $dX^0_1(s)/ds=L^0_1$,
where $L^0_1=1-1/4D_{LL}-1/2D_{L\un L}-1/4E_{LL}-1/2E_{L\un L}$, we can derive $\ds\f{1}{2}\le dt/ds\le 2$ under the assumptions
$|D|\leq 1/16$ and $|E|\leq 1/32$. Hence  integrating (2.30) from $s_1$ to $0$ together with (2.31) and Lemma 2.2 yields
\begin{align*}
&|\phi(t,x)|\le |\phi(X_1(s_1))||e^{G_1(s_1)}|+\ds\f{1}{2}\int^0_{s_1}\lt|e^{G_1(s)}
\rt|\lt|\lt((2L^i_1\p_i+\ds\f{\ell}{r})\phi\rt)(X_1(s))\rt|ds\\
&\le C|\phi(X_1(s_1))|+C\ds\int^0_{s_1}\lt|\lt((2L^i_1\p_i+\ds\f{\ell}{r})\phi\rt)(X_1(s))\rt|ds\\
&\le C\sup\limits_{\tau_1\leq\tau\le t}\ds\sum_{0\le I\leq 1}\|Z^I u(\tau,.)\|_{L^\infty}+C\int^t_{\tau_1}\bigl((1+\tau)\|\wt\sq_g u(\tau,.)\|_{L^\infty(H_{\tau})}\\
&\qquad +\ds\sum_{0\le I\leq 2}(1+\tau)^{-1}\|Z^I u(\tau,.)\|_{L^\infty(H_{\tau})}\bigr)d\tau,\tag{2.32}
\end{align*}
here we have used the following fact
$$|\phi(X_1(s_1))|\le C|(r-t){\un L}u (X_1(s_1))|\le C\ds\sum_{0\le I\le 1}|Z^I u(X_1(s_1))|\le C\sup\limits_{\tau_1\leq\tau\le t}\ds\sum_{0\le I\leq 1}\|Z^I u(\tau,.)\|_{L^\infty}$$
due to $X_1(s_1)\in \p H$ (which means $r=t/2 $ or  $3t/2$, and then  $r=t-r$  or  $3(r-t)$
holds), $\p_q=-\f12{\un L}$  and  ${\un L}=\ds\f{S-\o^\al\G_{0\al}}{r-t}$.

If $\tau_1<1$, then we know that there exists $\bar s$ with $s_1<\bar s\le 0$ such that  $X^0_1(\bar s)=1$
since  $t$  is decreasing along the backward integral curve $X_1(s)$  and  $s$ is an
increasing function  of  $t$.
Thus, as in (2.32), integrating (2.30) from $\bar s$ to $0$ yields
\begin{align*}
&|\phi(t,x)|\le |\phi(X_1(\bar s))||e^{G_1(\bar s)}|+\ds\f{1}{2}\int^0_{\bar s}\lt|e^{G_1(s)}\rt|\lt|\lt((2L^i_1\p_i+\ds\f{\ell}{r})\phi\rt)(X_1(s))\rt|ds\\
&\le C|\phi(X_1(\bar s))|+\ds\int^0_{\bar s}\lt|\lt((2L^i_1\p_i+\ds\f{\ell}{r})\phi\rt)(X_1(s))\rt|ds\\
&\le C\sup\limits_{\tau_1\leq\tau\le t}\ds\sum_{0\le I\leq 1}\|Z^I u(\tau,.)\|_{L^\infty}+C\int^t_{1}\bigl((1+\tau)\|\wt\sq_g u(\tau,.)\|_{L^\infty(H_\tau)}\\
&\qquad +\ds\sum_{0\le I\leq 2}(1+\tau)^{-1}\|Z^I u(\tau,.)\|_{L^\infty(H_\tau)}\bigr)d\tau\\
&\le C\sup\limits_{\tau_1\leq\tau\le t}\ds\sum_{0\le I\leq 1}\|Z^I u(\tau,.)\|_{L^\infty}+C\int^t_{\tau_1}\bigl((1+\tau)\|\wt\sq_g u(\tau,.)\|_{L^\infty(H_\tau)}\\
&\qquad +\ds\sum_{0\le I\leq 2}(1+\tau)^{-1}\|Z^I u(\tau,.)\|_{L^\infty(H_\tau)}\bigr)d\tau,\tag{2.33}
\end{align*}
here we have applied the facts of $r\le 1+X^0_1(\bar s)=2$ and
$$|\phi (X_1(\bar s))|\le 2|\p_q u|\le C \ds\sum_{0\le I\le 1}|Z^I u(X_1(\bar s))|
\le C \sup\limits_{\tau_1\leq\tau\le t}\ds\sum_{0\le I\leq 1}\|Z^I u(\tau,.)\|_{L^\infty}.$$
Therefore, combining (2.29) with (2.32)-(2.33) yields (2.26).

To prove (2.28), one can follow similarly from integrating $\ds\f{d}{ds}\phi(X_2(s))=(L^i_2\p_i\phi)(X_2(s))$ and
applying (2.16) in Lemma 2.2. \qquad\qquad\qquad\qquad\qquad\qquad\qquad$\square$

\vskip 0.3 true cm
As in [20], let $\r=\r(t,x)$ be constant along the integral curves of the  vector field $L_2=L^i_2\p_i$ close to the light cone and
equal to $r-t$ outside a neighborhood of the
forward light cone so that the usual phase $r-t$ can be replaced by $\rho(t,x)$ for the equation (1.3)
on a curved background. Namely, $\r=\r(t,x)$ satisfies
$$L^i_2\p_i\r=0 \;\;\text{when}\quad |t-r|\le t/2 , \;\;\;\r=r-t \;\;\text{when}\quad |t-r|\ge t/2.\eqno{(2.34)}$$
We notice that $(\tau_2, \bar x)$ is the first intersection point of the backward integral curve with $|r-t|=t/2$, then
by the definition of $\rho$, we have
$$|\rho(t,x)|=|\rho(\tau_2,\bar x)|=\tau_2/2\le t/2\qquad\text{when $|t-r|\le t/2$}.$$
In addition, we will take $\rho=\rho(q,p,\o)$ as a function of $q=r-t, p=r+t$ and $\o=x/|x|$.
As can be shown below, $0<\p_q\r=\r_q<\infty$ holds, then  $q$ can be also considered an invertible function of $\r$ for fixed $(p, \o)$
and $\p_q=\r_q\p_\r$. We also note that
$$ [L^i_2\p_i , \p_q]=-\ds\f{\p_q D_{LL}}{2}\p_q-\ds\f{\p_q E_{LL}}{2}\p_q \qquad\text{and}\;\;[L^i_2\p_i, \p_\r]=0.\eqno{(2.35)}$$
Here we point out that by comparison with (5.5) in [20], there is an extra troublesome term $\ds\f{\p_q E_{LL}}{2}\p_q$ in  $[L^i_2\p_i, \p_q]$
of (2.35), which should be  specially paid attention since $\p_q E_{LL}$ contains the second order derivatives $\p^2u$
but $\p_q D_{LL}$ only contains the first order derivatives $\p u$. With respect to the technical treatments on the coefficient
$\p_q E_{LL}$, one can refer to $\S 3$ (see (3.9)-(3.11)) below.
Next  we list an equivalence relation between $\rho$ and $q$, and the estimate
on $\p_\r\p_q\r$ which are completely
analogous to the ones in Proposition 5.1 and Lemma 5.2 of [20] respectively.

{\bf Lemma 2.4.} {\it Let $\r(t,x)$ be defined as in (2.34), and assume that $D_{LL}$ and $E_{LL}$ satisfy
$$|\p_q D_{LL}|+|\p_q E_{LL}|\le \ds\f{C\ve}{1+t}\ds\f{1}{(1+|\r|)^\nu} \eqno{(2.36)}$$
and
$$|D_{LL}|+|E_{LL}|\le C\ve \ds\f{1+|q|}{1+t}\eqno{(2.37)}$$
for some $\nu \ge 0$. Then
$$\lt(\ds\f{1+t}{1+|\r|}\rt)^{-C\ve V(\r)}\le \r_q\le \lt(\ds\f{1+t}{1+|\r|}\rt)^{C\ve V(\r)}
\quad\text{with $V(\r)=(1+|\r|)^{-\nu}$}\eqno{(2.38)}$$
and
$$ \lt(\ds\f{1+t}{1+|\r|}\rt)^{-C\ve}\le \ds\f{1+|q|}{1+|\r|}\le \lt(\ds\f{1+t}{1+|\r|}\rt)^{C\ve}.\eqno{(2.39)}$$
In addition, if one further assumes
$$|\p^2_q D_{LL}|+|\p^2_q E_{LL}|\le \ds\f{C\r_q\ve}{1+t}\ds\f{1}{(1+|\r|)^{1+\nu}}\eqno{(2.40)}$$
for some $\nu\ge 0$, then
$$|\p_{\r}\p_q \r|\le \ds\f{C\r_q\ve}{(1+|\r|)^{1+\nu}}\ln|\ds\f{1+t}{1+|\r|}|.\eqno{(2.41)}$$}

\vskip 0.2 true cm
{\bf Remark 2.2.} {\it  By $E^{ij}=\ds\sum^3_{k=1}e^{ij}_k\p_k u$, it seems that the decay property on the time $t$
of $\p_q E_{LL}$ in (2.36) should coincide with that of the second order derivative $\p^2u$, namely, $|\p_q E_{LL}|\le C|\p^2u|\le C\ve (1+t)^{-1+C\ve}$ (see (3.61) of Proposition 3.9 below).  However,
thanks to the null condition property of $e^{ij}_k\p_k u\p_{ij}u$, we can show that $\p_q E_{LL}$
will admit better decay rate of $(1+t)^{-1}$ (see (3.9) in $\S 3$). This is one of the key points
that we can show the global existence of the solution to (1.3).}

\vskip 0.2 true cm
{\bf Remark 2.3.} {\it Here we point out that we only need the assumption (2.36) to derive (2.38). And the assumption (2.37) suffices to get (2.39).  }

\vskip 0.2 true cm

{\bf Proof.} The proof procedures are completely similar to those in Proposition 5.1 and Lemma 5.2 of [20] under the assumptions
(2.36)-(2.37) and (2.40), then we omit the details here.
\qquad\qquad  $\square$

\vskip 0.2 true cm
At the end of this section, we will give two useful  inequalities
when the related null condition holds.

{\bf Lemma 2.5.} {\it Assume $\ds\sum^3_{k=0}e^{ij}_k\o_i\o_j\o_k
\equiv0$ (i.e., the null condition) holds,
then for smooth functions $u, v$ supported in $|x|\le 1+t$, we have
$$ |\ds\sum^3_{k=0}e^{ij}_k \p_k u\p^2_{ij}v|\le C (|\bar Z u||\p^2 v|+|\p u||\bar Z \p v|)\eqno{(2.42)}$$
and
$$|\ds\sum^3_{k=0} e^{ij}_k\p_k u\p^2_{ij} v|\le C(1+t)^{-1}(|Zu||\p^2 v|+|\p u||Z\p v|),\eqno{(2.43)}$$
where $\bar Z=\lt\{\p_1+\o_1\p_t,\p_2+\o_2\p_t,\p_3+\o_3\p_t\rt\}$.}

{\bf Proof.} By the null condition $\ds\sum^3_{k=0}e^{ij}_k\o_i\o_j\o_k\equiv0$, then
\begin{align*}
&\ds\sum^3_{k=0}e^{ij}_k\p_k u \p^2_{ij}v\\
&=\ds\sum^3_{k=0}\biggl(e^{ij}_k(\p_k u+\o_k\p_t u)\p^2_{ij}v-e^{ij}_k\o_k\p_t u(\p_i+\o_i\p_t)\p_j v
+e^{ij}_k\o_k\o_i\p_t u(\p_j+\o_j\p_t)\p_t v\\
&\qquad  -
e^{ij}_k\o_k\o_i\o_j\p_t u\p^2_t v\biggr)\\
&=\ds\sum^3_{k=0}\biggl(e^{ij}_k(\p_k u+\o_k\p_t u)\p^2_{ij}v-e^{ij}_k\o_k\p_t u(\p_i+\o_i\p_t)\p_j v
+e^{ij}_k\o_k\o_i\p_t u(\p_j+\o_j\p_t)\p_t v\biggr), \end{align*}
which derives (2.42).

In addition, by $\p_k u+\o_k \p_t u=t^{-1}\G_{0k}u-\o_k t^{-1}(r-t)\p_t u$ ($k=1, 2, 3$) and (2.6), we then have
$$|\ds\sum^3_{k=0} e^{ij}_k\p_k u\p^2_{ij} v|\le C(1+t)^{-1}(|Zu||\p^2 v|+|\p u||Z\p v|).$$
Therefore, the proof of Lemma 2.5 is completed.\qquad\qquad \qquad \qquad \qquad  $\square$

\vskip 0.5 true cm
\centerline{\bf $\S 3$. The sharp decay estimate of the solution $u$ to (1.3)}
\vskip 0.4 true cm
As in [20], we assume that the solution $u$ of (1.3) admits the following weak decay estimate
$$ |Z^I u|\le M\ve(1+t)^{-\nu},\;\;\;|I|\le N-3 ,\;\;M\ve\le 1\eqno{(3.1)}$$
for some large $N\ge 8$,  $\ds\f12<\nu<1$ and a positive constant $M$.  From this, we manege to derive the strong
decay estimate of $u$ and further
obtain the more precise energy estimates in $\S 5$ and $\S 6$. Below we denote $C>0$ by a generic
constant depending only on $M$.

By the finite propagation speed property
for the wave equation (1.3), then
$$u(t,x)\equiv 0\qquad \text{for $r\ge 1+t$}. \eqno{(3.2)}$$
In addition, by scaling we may assume from now on for $0<c_0<<1$
$$|D|+|E|\le \ds c_0(|u|+|\p u|).\eqno{(3.3)}$$

\vskip 0.2 true cm
First, we establish the strong decay estimates of $u$ and $\p u$ than the ones in (3.1).

{\bf Lemma 3.1.} {\it Assume that the solution $u$ to (1.3) satisfies (3.1)-(3.3).
Then
\begin{equation*}
|u|\le\left\{
\begin{aligned}
&C \ve (1+t)^{-1}(1+|q|), \\
&C \ve (1+t)^{-1+C \ve}(1+|\rho|)^{1-\nu-C\ve}
\end{aligned}
\right.\tag{3.4}
\end{equation*}
and
$$|\p u|\le C\ve (1+t)^{-1}(1+|\r|)^{-\nu},\eqno{(3.5)}$$
where $q=r-t$, and the function $\rho$ has been defined in (2.34).}

\vskip 0.2 true cm
{\bf Proof.} At first, we prove the first estimate in (3.4). By (3.1) and (3.3), then for any $T>0$
and small $\ve$
\begin{align*}
&\int^T_0 (\|D(t,.)\|_{L^\infty(H_t)}+\|E(t,.)\|_{L^\infty(H_t)})\ds\f{dt}{1+t}\le c_0\int^{\infty}_0 \ds\f{(|u|+|\p u|)}{1+t}dt\\
&\le
c_0 M\ve\int^{\infty}_0 (1+t)^{-1-\nu} dt \leq 1,
\end{align*}
where  $H_t=\{x\in\Bbb R^3: \ds\f{t}{2}<|x|<\ds\f{3t}{2}\}$.
Together with (2.26), (3.1) and $\wt\sq_g u=0$, this yields
$$(1+t+r)|\p u|\le C\ve$$
and
$$|\p u|\le C\ve (1+t+r)^{-1}\le C\ve (1+t)^{-1}.\eqno{(3.6)}$$
Thus, the first estimate of (3.4) follows from integrating (3.6) from $r=1+t$ where $u=0$,
$$|u(t,x)|=|u(t,r,\o)|=|-\int^{1+t}_r\p_{r'} u (t,r',\o)dr'|\le \int^{1+t}_r C\ve(1+t)^{-1}dr'
\le C\ve(1+t)^{-1}(1+|q|).\eqno{(3.7)}$$

Next, we show (3.5). In fact, by (3.6) and (3.7), we have
$$|\p D|+|E|+(1+|q|)^{-1}|D|\le \ds\f{C\ve}{1+t}$$
or
$$|\p D|+|E|+(1+|q|)^{-1}|D|\le \ds\f{C\ve}{1+t+r}.$$
Therefore, by choosing $a=\mu=0$ in Lemma 4.2 of [20] (note that we have established Lemma 2.3 in $\S 2$,
then the conclusion analogous to Lemma 4.2 of [20] follows by an easy computation), we can obtain
\begin{align*}
&(1+t+r)(1+|\r(t,x))^{\nu}|\p u(t,x)|\leq C \sup\limits_{0\leq\tau\le t}(1+\tau)^\nu \ds\sum_{0\le I\le 2}\|Z^I u(\tau,.)\|_{L^\infty}\\
&\qquad\quad +C\int^t_{0}(1+\tau)\|(1+|\r(\tau,.)|)^{\nu}\wt\sq_g u(\tau,.)\|_{L^\infty(H_\tau)} d\tau\\
&\qquad =C \sup\limits_{0\leq\tau\le t}(1+\tau)^\nu \ds\sum_{0\le I\le 2}\|Z^I u(\tau,.)\|_{L^\infty}.
\end{align*}
This, together with (3.1), yields (3.5).

Finally, we show the second estimate in (3.4).
It is noted that we have by (3.5)
$$ |\p_q D_{LL}|+|E_{LL}|\le \ds\f{C\ve}{1+t}\ds\f{1}{(1+|\r|)^\nu}.\eqno{(3.8)}$$
To derive the second estimate in (3.4), we will apply Lemma 2.4 and (3.5)
to derive such a kind of estimate $|\p u|\le C\ve(1+t)^{-1+C\ve}(1+|q|)^{-\nu}$
so that the estimate of $u$ can be obtained by integrating from $r=t+1$ where $u=0$.
To this end, we require to verify the assumptions (2.36)-(2.37) of  Lemma 2.4. In fact, by (3.7)-(3.8), we only need to verify
$$|\p_q E_{LL}|\le \ds\f{C\ve}{1+t}\ds\f{1}{(1+|\r|)^\nu}.\eqno{(3.9)}$$
Indeed, by the null condition  $\ds\sum^3_{k=0}e^{ij}_k\o_i\o_j\o_k\equiv 0$ and  $L_i=\o_i$
(due to $L_0=-L^0=-1$ and $L_{\al}=\o_{\al}$), we have
\begin{align*}
\p_q E_{LL}
&=-\ds\f{1}{2}{\un L}^m\p_m(E^{ij}L_i L_j)\\
&=-\ds\f{1}{2}{\un L}^m\p_m(\ds\sum^3_{k=0}e^{ij}_k\p_k u L_i L_j)\\
&=-\ds\f{1}{2}{\un L}^m \ds\sum^3_{k=0}e^{ij}_k(\p_m\p_k u)L_i L_j\\
&=-\ds\f{1}{2}{\un L}^m e^{ij}_0(\p_t\p_m u+\o_m\p^2_t u)L_i L_j-\ds\f{1}{2}{\un L}^m \ds\sum^3_{\al=1}e^{ij}_{\al}(\p_{\al}\p_m u
-\o_{\al}\o_m \p^2_t u)L_i L_j\\
&\;\;\;\;-\ds\f{1}{2}{\un L}^m\o_m (-e^{ij}_0 L_i L_j+e^{ij}_{\al} \o_{\al} L_i L_j)\p^2_t u \\
&=-\ds\f{1}{2}{\un L}^m e^{ij}_0(\p_t\p_m u+\o_m\p^2_t u)L_i L_j-\ds\f{1}{2}{\un L}^m \ds\sum^3_{\al=1} e^{ij}_{\al}(\p_{\al}\p_m u
-\o_{\al}\o_m \p^2_t u)L_i L_j,\tag{3.10}
\end{align*}
here we have used the crucial null condition to derive $-e^{ij}_0 L_i L_j+e^{ij}_{\al} \o_{\al} L_i L_j\equiv 0$.

In addition, one also has  for $\al =1,2,3$
$$\p_t\p_{\al}u+\o_{\al}\p^2_t u=t^{-1}\G_{0\al}\p_t u -t^{-1}\o_{\al}q\p^2_t u,$$
$$\p_1\p_{\al}u-\o_1\o_{\al}\p^2_t u=t^{-1}(\G_{01}\p_{\al}u -\o_1\G_{0\al}\p_t u)-t^{-1}\o_1q(\p_t\p_{\al}u
-\o_{\al}\p^2_t u),$$
and $\p_2\p_{\al}$, $\p_3\p_{\al}$  proceed similarly. According to  this and (3.10), together with (3.1)
and the facts of $|\p^2 u|\le C (1+|q|)^{-1}|Z \p u|$  and  $|\r|\le Ct$, we arrive at
$$ |\p_q E_{LL}|\le C (1+t)^{-1}(|Z\p u|+|q| |\p^2 u|)\le C (1+t)^{-1}|Z\p u|\le C\ve (1+t)^{-1}(1+|\r|)^{-\nu},\eqno{(3.11)}$$
that is, (3.9) is shown.

By (3.8) and (3.9), we hence conclude that  the assumptions (2.36)-(2.37) of Lemma 2.4 hold.  Then by (2.39)
in Lemma 2.4, one has
$$|\p u|\le C\ve (1+t)^{-1}(1+|\r|)^{-\nu}\le C\ve (1+t)^{-1+C\ve}(1+|q|)^{-\nu}.$$
From this and $\ds\int^{1+t}_r (1+|r'-t|)^{-\nu}dr'\le C(1+|q|)^{1-\nu}$, we have
\begin{align*}
&|u(t,x)|=|-\int^{1+t}_r \p_{r'} u (t,r',\o)dr'|\le C\ve (1+t)^{-1+C\ve}\int^{1+t}_r (1+|r'-t|)^{-\nu}dr'\\
&\le C\ve (1+t)^{-1+C\ve}(1+|q|)^{1-\nu}.\tag{3.12}
\end{align*}
Using (2.39) again for (3.12), we then get the second estimate in (3.4).  \qquad  \qquad  \qquad
 \qquad \qquad $\square$

Next we derive the decay estimate of $\p^2 u$. Before doing this, as in $\S 6$ of [20], we require to establish
some Lemmas so that $L^i_2\p_i(r\p^2_q u)$ and $L^i_2\p_i(r\p_\r\p_q u)$ can be suitably approximated by
$r\p_q\wt\sq_g u$ and $r\p_\r\wt\sq_g u$ respectively (in fact, $r\p_q\wt\sq_g u=0$ and $r\p_\r\wt\sq_g u=0$
hold by the equation (1.3)). If so, integrating along the integral curve
of $L_2^i\p_i$, one can obtain the decay of $\p^2 u$. In this process, we should specially pay attention to
the terms $E$ and $\p E$ which include $\p u$ and $\p^2 u$ respectively.

{\bf Lemma 3.2.} {\it Assume that for some constant $a\ge 0$
$$|D|\le \ds\f{1+|t-r|}{1+t+r}\lt|\ds\f{1+t+r}{1+|t-r|}\rt|^a\eqno{(3.13)}$$
and
$$|\p D|+|E|\leq \ds\f{1}{1+t+r}\lt|\ds\f{1+t+r}{1+|t-r|}\rt|^a.\eqno{(3.14)}$$
Then
\begin{align*}
|2L^i_2\p_i(r\p^2_q u)+r(\p_q D_{LL}+\p_q E_{LL})\p^2_q u-r\p_q\wt\sq_g u|&\leq
\ds\f{C}{1+t+r}\lt|\ds\f{1+t+r}{1+|t-r|}\rt|^a \ds\sum_{0\le I\leq 2}|\p Z^I u|\\
&\quad +C|\p E||\p Z u|.\tag{3.15}
\end{align*}}

\vskip 0.2 true cm
{\bf Remark 3.1.} {\it Compared with (3.15) of Lemma 3.3 in [20], here an extra term $ |\p E||\p Z u|$
appears in the right hand side of (3.15). This term cannot be controlled
by $\ds\f{C}{1+t+r}\lt|\ds\f{1+t+r}{1+|t-r|}\rt|^a$ $\times \ds\sum_{0\le I\leq 2}|\p Z^I u|$ directly
since we have no estimate on $\p E$ so far. In addition,  the coefficient
of $\p^2_q u$ in the left hand side of (3.15) is also different from that in (3.15) of [20] due to the appearance
of $\p_q E_{LL}$.}
\vskip 0.2 true cm
{\bf Proof.} It follows from (2.4) and a direct computation (or see [20, Lemma 2.2]) that for $f^{ij}=\p_m(D^{ij}+E^{ij})$
$$|f^{ij}\p_i\p_j u-f_{LL}\p^2_q u| \leq C |\bar\p\p u|,\eqno{(3.16)}$$
where $C=\sup\limits_{i,j}|f^{ij}|$ and
$\bar\p=\{L, S_1, S_2\}$.

In addition, one has
\begin{align*}
f_{LL}&=f^{ij}L_i L_j=\p_m(D^{ij}+E^{ij})L_i L_j\\
&=\p_m D_{LL}-D^{ij}(\p_m L_i)L_j-D^{ij}L_i(\p_m L_j)+\p_m E_{LL}-E^{ij}(\p_m L_i)L_j-E^{ij}L_i(\p_m L_j).\tag{3.17}
\end{align*}
Note that $|\p_m L_i|\le C/r$  and $|\p^2_q u|\le C (1+|t-r|)^{-1}|\p Z u|$ hold. This,
together with (2.5), (3.16)-(3.17) and the assumptions (3.13)-(3.14), yields
\begin{align*}
& |\p_m(D^{ij}+E^{ij})\p_i\p_j u-(\p_m D_{LL}+\p_m E_{LL})\p^2_q u|\\
\le \;&C(|\p D|+|\p E|)|\bar\p \p u|+|\lt(D^{ij}(\p_m L_i)L_j+D^{ij}L_i(\p_m L_j)
+E^{ij}(\p_m L_i)L_j+E^{ij}L_i(\p_m L_j)\rt)\p^2_q u|\\
\le\; &\ds\f{C}{r}\ds\f{1}{1+t+r}\lt|\ds\f{1+t+r}{1+|t-r|}\rt|^a\ds\sum_{0\le I\le 1}|\p Z^I u|
+\ds\f{C}{1+t+r}|\p E||\p Z u|.\tag{3.18}
\end{align*}
On the other hand, if we use  $\p_m u$ in (2.16) instead of $u$,  then
$$ |2L^i_2\p_i(r\p_q \p_m u)-r\wt\sq_g \p_m u|\leq \ds\f{C}{1+t+r}\lt|\ds\f{1+t+r}{1+|t-r|}\rt|^a \ds\sum_{0\le I\leq 2}|Z^I \p_m u|.\eqno{(3.19)}$$
Since
$$\p_m \wt\sq_g u=\wt\sq_g\p_m u+(\p_m D^{ij})\p_i\p_j u+(\p_m E^{ij})\p_i\p_j u,\eqno{(3.20)}$$
by (3.19)-(3.20) and (3.18) we have
\begin{align*}
&|2L^i_2\p_i(r\p_q \p_m u)+r(\p_m D_{LL}+\p_m E_{LL})\p^2_q u-r\p_m\wt\sq_g u|\\
\le\; &|2L^i_2\p_i(r\p_q \p_m u)-r\wt\sq_g \p_m u|+r|-\p_m(D^{ij}+E^{ij})\p_i\p_j u+(\p_m D_{LL}+\p_m E_{LL})\p^2_q u|\\
\le \;&\ds\f{C}{1+t+r}\lt|\ds\f{1+t+r}{1+|t-r|}\rt|^a \ds\sum_{0\le I\leq 2}|Z^I \p_m u|+C|\p E||\p Z u|.\tag{3.21}
\end{align*}
Using $\p_q=-\ds\f{1}{2}{\un L}=-\ds\f{1}{2}{\un L} ^m\p_m$ in (3.21), then we complete the proof of (3.15).\qquad \qquad \qquad \qquad
$\square$

{\bf Lemma 3.3.} {\it Assume
$$|\p D|+|E|+(1+|q|)^{-1}|D|\le \ds\f{C\ve}{1+t+r},\eqno{(3.22)}$$
then
$$|2L^i_2\p_i(r\p_\r\p_q u)-r\p_\r\wt\sq_g u|\leq \ds\f{C\r^{-1}_q}{1+t+r}| \ds\sum_{0\le I\leq 2}|\p Z^I u|+C\r^{-1}_q|\p E||\p Z u|\eqno{(3.23)}$$
and
$$|2L^i_2\p_i(r\p_\r\phi)-r\r^{-1}_q\wt\sq_g \phi|\leq c_1\ve|\p_\r\phi|+\ds\f{C\r^{-1}_q}{1+t+r}| \ds\sum_{0\le I\leq 2}|\p Z^I\phi|.\eqno{(3.24)}$$}
\vskip 0.2 true cm
{\bf Remark 3.2.} {\it Similar to Remark 3.1, by comparison with (6.11) of Lemma 6.2 in [20], here an extra term $\r^{-1}_q|\p E||\p Z\phi|$
also appears in the right hand side of (3.23).}
\vskip 0.2 true cm

{\bf Proof.} Due to $L^i_2\p_i \r=0$, $2 L^i_2\p_i\p_q\r=-\p_q D_{LL}\p_q \r-\p_q E_{LL}\p_q \r$ and  $2 L^i_2\p_i \r^{-1}_q=\r^{-1}_q\p_q D_{LL}+\r^{-1}_q\p_q E_{LL}$,
it follows that
$$2L^i_2\p_i(r\r^{-1}_q\p^2_q u)=\r^{-1}_q\lt[2L^i_2\p_i(r\p^2_q u)+(\p_q D_{LL}+\p_q E_{LL})r\p^2_q u\rt]\eqno{(3.25)}$$
and further
$$ 2L^i_2\p_i(r\p^2_q u)=2\r_q L^i_2\p_i(r\r^{-1}_q\p^2_q u)-r(\p_q D_{LL}+\p_q E_{LL})\p^2_q u.\eqno{(3.26)}$$
Applying Lemma 3.2 with $a=0$ for (3.26) yields
$$|2\r_q L^i_2\p_i(r\r^{-1}_q\p^2_q u)-r\p_q\wt\sq_g u|\le \ds\f{C}{1+t+r} \ds\sum_{0\le I\leq 2}|\p Z^I u|+C|\p E||\p Z u|,$$
which derives (3.23).

(3.24) follows from (2.16) with $a=0$ and $|\p_q E_{LL}|\le C\ve (1+t)^{-1}$ directly.\qquad\qquad $\square$
\vskip 0.3 true cm

Based on Lemma 3.2 and Lemma 3.3, we now establish the strong decay estimate of second order derivatives $\p^2 u$.

{\bf Lemma 3.4.} {\it Assume that the solution $u$ to (1.3) satisfies (3.1)-(3.3).
Then for small $\ve>0$
$$|\p^2 u|\le C\ve (1+t)^{-1}(1+|\r|)^{-1-\nu}\r_q.\eqno{(3.27)}$$
}

{\bf Proof.} When $(t,x)\not\in H$, $|q|\ge t/2$ and $\r=q$ and further $\ds\f{1+t}{1+|\r|}=O(1)$. Then it follows from this, (2.6) and (3.1)
that
\begin{align*}
|\p^2 u|&\le C(1+|q|)^{-2}\ds\sum_{0\le I\le 2}|Z^I u|\\
&\le C\ve(1+|\r|)^{-2}(1+t)^{-\nu}\\
&\le C\ve (1+|\r|)^{-1-\nu}\lt(\ds\f{1+t}{1+|\r|}\rt)^{1-\nu}(1+t)^{-1}\\
&\le C\ve (1+t)^{-1}(1+|\r|)^{-1-\nu},
\end{align*}
which implies (3.27) holds.

We now consider the case of $(t,x)\in H$.
Notice that
$\p^2$ can be expressed as the combinations of $\p^2_q $, $\bar\p^2 $ and $\bar\p \p_q$ with bounded coefficients, and thus
$$ |\p^2 u|\le C (|\p^2_qu|+|\bar\p^2 u|+|\bar\p \p_q u|).\eqno{(3.28)}$$

We now analyze each term in the right hand side of (3.28).
To obtain the estimate of $\p^2_q u$, we first study $L^i_2\p_i(r\rho^{-1}_q\p^2_q u)$. For this aim,
applying (3.23), Lemma 2.4 and (3.1) yields
\begin{align*}
|2L^i_2\p_i(r\rho^{-1}_q\p^2_q u)|
&\le \ds\f{C\rho^{-1}_q}{1+t}\ds\sum_{0\le I\le 2}|\p Z^I u|+C\rho^{-1}_q|\p E||\p Z u|\\
&\le \ds\f{C\rho^{-1}_q}{(1+t)(1+|q|)}\ds\sum_{0\le I\le 3}|Z^I u|+\ds\f{C\rho^{-1}_q}{(1+|q|)^2}|Z E||Z^2 u|\\
&\le \ds\f{C\ve}{(1+t)^{1+\nu-C\ve}(1+|\rho|)^{1+C\ve}}+\ds\f{C\ve}{(1+t)^{2\nu-C\ve}(1+|\rho|)^{2+C\ve}}.\tag{3.29}
\end{align*}
Since $\nu >\ds\f{1}{2}+C\ve$ holds for small $\ve$,  integrating (3.29) from $r=t/2$
(at this place, $t\sim|\r|$  and  $|r\r^{-1}_q\p^2_q u|\le C\ve(1+|\r|)^{-\nu-1}$) yields
$$|r\r^{-1}_q\p^2_q u|\le C\ve(1+|\r|)^{-\nu-1}$$
and further
$$|\p^2_q u|\le C\ve (1+t)^{-1}(1+|\r|)^{-\nu-1}\r_q.\eqno{(3.30)}$$

Next we estimate $\bar\p^2 u$. Without loss of generality, we assume $t>1$. In this case, $r>\ds\f{t}{2}>\ds\f{1+t}{4}$
holds for $(t,x)\in H$. Therefore, by (2.8) and (3.1), we have
\begin{align*}
|\bar\p^2 u|&\le \ds\f{C}{r}\ds\sum_{0\le I\le 2}\ds\f{|Z^I u|}{1+t+r}\\
&\le C\ve(1+t)^{-1}(1+t)^{-1-\nu}\\
&\le C\ve(1+t)^{-1}(1+|\rho|)^{-1-\nu}\rho_q.\tag{3.31}
\end{align*}
Similarly, we have
$$|\bar\p\p_q|\le C\ve(1+t)^{-1}(1+|\rho|)^{-1-\nu}\rho_q.$$
This, together with (3.30)-(3.31) and (3.28), yields (3.27).
\qquad\qquad\qquad\qquad\qquad $\square$

Next we establish the strong decay estimate on $Zu$. To this end, we require
to calculate the commutators of vector fields $Z$  with
$\wt\sq_g=\sq+D^{ij}\p_i\p_j+\ds\sum^3_{k=0}e^{ij}_k\p_k u\p_i\p_j$.
Below $G(u,v)$ denotes the various bilinear form analogous to\; $e^{ij}_k\p_k u\p^2_{ij} v$, which satisfies the null condition.
It follows from Lemma 6.6.5 in [11] and a direct computation that
\begin{align*}
&Z\wt\sq_g \phi= Z\sq\phi+Z (D^{ij}\p_i\p_j\phi)+Z G(u,\phi)\\
&=\sq Z\phi-C_Z\sq\phi+D^{ij}\p_i\p_j Z\phi+(Z D^{ij})\p_i\p_j\phi+2D^{ij}C^l_{Zi}\p_l\p_j\phi+G(Zu,\phi)\\
&\qquad \qquad \qquad\qquad \qquad\qquad \qquad\qquad \qquad\qquad \qquad +G(u,Z\phi)+G(u,\phi)\\
&=\wt\sq_g Z\phi-C_Z\wt\sq_g\phi+C_Z D^{ij}\p_i\p_j\phi+2D^{ij}C^l_{Zi}\p_l\p_j\phi+(Z D^{ij})\p_i\p_j\phi+G(Zu,\phi)+G(u,\phi),
\end{align*}
where $C_Z$ and $C^l_{Zi}$ are some suitable constants. Set $\hat Z=Z+C_Z$, then we have
$$\wt\sq_g Z\phi=\hat Z\wt\sq_g\phi-(ZD^{ij}+2C^i_{Zl}D^{lj}+C_Z D^{ij})\p_i\p_j\phi+G(Zu,\phi)+G(u,\phi)\eqno{(3.32)}$$
and further
$$\wt\sq_g Z^I \phi=\hat Z ^I\wt\sq_g\phi+\ds\sum_{J+K\le I,K<I}C^{\;I \;\;lm}_{JK ij}(Z^J D^{ij})\p_l\p_m Z^K \phi+\ds\sum_{J+K\le I,K<I}G(Z^J u,Z^K \phi),\eqno{(3.33)}$$
where $C^{\;I \;\;lm}_{JK ij}$ are constants.

More generally, we can have
$$\wt\sq_g \p^kZ^I \phi=\p^k\hat Z^I\wt\sq_g\phi+\ds\sum_{J+K\le I, s+n=k, n+K<k+I}C^{\;I \;\;lm}_{JK ij}(\p^s Z^J D^{ij})\p_l\p_m \p ^n Z^K \phi$$
$$+\ds\sum_{J+K\le I, s+n=k, n+K<k+I}G(\p^s Z^J u,\p^n Z^K \phi).\eqno{(3.34)}$$
Thus, by (3.34) and the equation (1.3), we have for $u$
$$|\wt\sq_g \p^k Z^I u|\le C\ds\sum_{J+K\le I,K<I}|Z^J u||\p^2\p^k Z^K u|+\ds\sum_{s+n=k,J+K\le I}|\p\p^s Z^J u|
|\p\p^n Z^K u|$$
$$+\ds\sum_{J+K\le I, s+n=k, n+K<k+I}|G(\p^s Z^J u,\p^n Z^K \phi)|.\eqno{(3.35)}$$

Before estimating $Zu$,  we require to cite the result in Lemma 6.3 of [20] for reader's convenience.
Here we point out that although the form of $L^i_2$ in our paper is somewhat different from that in [20],  by minor modification
on the proof of Lemma 6.3 in [20]
(still integrating along the integral curve of $L_2$ and applying Gronwall's inequality),  we have

{\bf Lemma 3.5.} {\it Assume for some constant $\nu>\ds\f12$
$$(1+|\r|)|\p_\r \phi|+|\phi|\le C\ve(1+|\r|)^{-\nu}\qquad \;\;when \;\;|t-r|=t/2\;\quad or\quad \;t+r\le 2$$
and
$$\phi=0\qquad when \quad r>1+t\quad and\quad t>0.$$
Assume also that for $(t,x)\in H$
$$|L^i_2\p_i(r\p_\r \phi)|\le C\ve \biggl(\ds\f{|\phi|}{1+|\r|}+|\p_\r \phi|\biggr)+\ds\f{C\ve^2}{(1+t)^{1-C\ve}(1+|\r|)^{\nu+C\ve}}+\ds\f{C\ve(1+|\r|)^{-C\ve}}{(1+t)^{1+\nu-C\ve}}.$$
Then
$$|\phi|(1+|\r|)^{-1}+|\p_\r \phi|\le C\ve(1+t)^{-1+C\ve}(1+|\r|)^{-\nu}.$$
}

With respect to the decay estimate of $Zu$, we have

{\bf Lemma 3.6.} {\it Assume that the solution $u$ to (1.3) satisfies (3.1)-(3.3).
Then
\begin{equation*}
|Zu|\le\left\{
\begin{aligned}
& C \ve (1+t)^{-1+C\ve}(1+|\rho|)^{1-\nu},\\
& \ve (1+t)^{-1}(|q|+(1+t)^{C\ve}).
\end{aligned}
\right.\tag{3.36}
\end{equation*}
}

{\bf Proof.} Due to $\wt\sq_g u=0$, we have by (3.32) and Lemma 2.5
$$|\wt\sq_g Z\phi|\le C(|Z u|+|u|)|\p^2 u|+C(1+t)^{-1}(\ds\sum_{1\le I\le 2}|Z^I u||\p^2 u|+\ds\sum_{0\le I\le 1}|\p Z^I u||Z\p u|).\eqno{(3.37)}$$
On the other hand, if we use $Z u$ instead of $\phi$ in (3.24) and combine (3.37), then
\begin{align*}
|L^i_2\p_i(r\p_\r Zu)|&\le C r(|Z u|+|u|)\ds\f{|\p^2 u|}{\r_q}+C\r^{-1}_q(\ds\sum_{1\le I\le 2}|Z^I u||\p^2 u|+\ds\sum_{0\le I\le 1}|\p Z^I u||Z\p u|)\\
&+C\ve|\p_\r Z u|+\ds\f{C\r^{-1}_q}{1+t}| \ds\sum_{0\le J\leq 3}|\p Z^J u|.
\end{align*}
Therefore it follows from (3.1), (3.4)-(3.5), (3.27) and (2.38)-(2.39) that
$$|L^i_2\p_i(r\p_\r Zu)|\le C\ve \biggl(\ds\f{|Z u|}{1+|\r|}+|\p_\r Zu|\biggr)+\ds\f{C\ve^2}{(1+t)^{1-C\ve}(1+|\r|)^{2\nu+C\ve}}+\ds\f{C\ve(1+|\r|)^{-C\ve}}{(1+t)^{1+\nu-C\ve}}.\eqno{(3.38)}$$
Note that  when  $|t-r|= t/2$
$$|Z u|+(1+|\r|)|\p_\r Zu|=|Zu|+(1+|q|)|\p Zu|\le C\ds\sum_{0\le I\le 2}|Z^I u|\le C\ve(1+t)^{-\nu}\le C\ve(1+|\r|)^{-\nu}.\eqno{(3.39)}$$
Therefore, by (3.38)-(3.39), one knows that all the assumptions of Lemma 3.5 are fulfilled, then it follows that
$$|Z u|(1+|\r|)^{-1}+|\p_\r Z u|\le C\ve(1+t)^{-1+C\ve}(1+|\r|)^{-\nu}.$$
Together with  $|\r|\le Ct$  and $(1+t)^{C\ve}(1+|q|)^{1-\nu}\le(1+t)^{C\ve/\nu}+(1+|q|)$, this yields
(3.36).\qquad \qquad $\square$
\vskip 0.2 true cm

Next, we establish the decay estimate of $\p^ku$ for $1\le k\le N-4$.

{\bf Lemma 3.7.} {\it Assume that the solution $u$ to (1.3) satisfies (3.1)-(3.3). For $1\le k \le N-4$, we have
$$ |\p^k u|\le
\ds\f{C\ve}{1+t}(1+|\r|)^{1-k-\nu}\lt(\ds\f{1+t}{1+|\r|}\rt)^{C\ve V(\rho)}, \eqno{(3.40)}$$
where $V(\rho)=(1+|\r|)^{-\nu}.$}

{\bf Proof.}   We will show (3.40) by induction method. If $k=1$, we have already proved  (3.40) by (3.5) of Lemma 3.1. Assume that
(3.40) holds for $k\le n\le N-5$, we will prove (3.40) for $k=n+1$. Taking $\p^n$ on two hand sides of $\wt\sq_g u=c^{ij}\p_i\p_j u+D^{ij}\p_i\p_j u+E^{ij}\p_i\p_j u$ yields
$$\wt\sq_g\p^n u=-\ds\sum_{m+k=n, m\ge1}(\p^m D^{ij})\p_i\p_j\p^k u-\ds\sum_{m+k=n, m\ge1}
(\p^m E^{ij})\p_i\p_j\p^ku.\eqno{(3.41)}$$
Substituting (3.5) and (3.40) with $k\le n$ into (3.41) yields
\begin{align*}
&|\wt\sq_g\p^n u|\\
&\le C|\p u|\ds\sum_{m=n+1}|\p^m u|+C\ds\sum_{2\le m\le n, 0\le k\le n-2, m+k=n}|\p^m u||\p^2\p^ku|
+\ds\sum_{m+k=n, m\ge1}|G(\p^m u,\p^k u)|\\
&\le \ds\f{C\ve}{(1+t)(1+|\r|)^\nu}\ds\sum_{m=n+1}|\p^m u|+\ds\f{C\ve^2}{(1+t)^2(1+|\r|)^{n+2\nu}}\lt(\ds\f{1+t}{1+|\r|}\rt)^{C\ve V(\r)}\\
&\qquad\qquad +\ds\sum_{m+k=n, m\ge1}|G(\p^m u,\p^k u)|.\tag{3.42}
\end{align*}
Note that
\begin{align*}
\ds\sum_{m+k=n, m\ge1}|G(\p^m u,\p^k u)|
&=\ds\sum_{k=n-1}|G(\p u,\p^k u)|+|G(\p^n u,u)|\\
&\quad +\ds\sum_{m+k=n, 2\le m\le n-1, 1\le k\le n-2}|G(\p^m u,\p^k u)|.\tag{3.43}
\end{align*}
By Lemma 2.5, (3.1), (2.6), (2.39) and the fact of $|\r|\le Ct$, we have
\begin{align*}
&\ds\sum_{k=n-1}|G(\p u,\p^k u)|+|G(\p^n u,u)|\\
\le\;& C(1+t)^{-1}(|Z\p u||\p^{n+1}u|+|\p^2 u||Z\p^n u|)\\
\le \;&\ds\f{C\ve}{(1+t)(1+|\r|)^\nu}\ds\sum_{m=n+1}|\p^m u|
+\ds\f{C\ve^2}{(1+t)^2(1+|\r|)^{n+2\nu}}\lt(\ds\f{1+t}{1+|\r|}\rt)^{C\ve V(\r)}.\tag{3.44}
\end{align*}
In addition, it  follows from (3.40)  with  $k\le n$ that
$$\ds\sum_{m+k=n, 2\le m\le n-1, 1\le k\le n-2}|G(\p^m u,\p^k u)|\le \ds\f{C\ve^2}{(1+t)^2(1+|\r|)^{n+2\nu}}\lt(\ds\f{1+t}{1+|\r|}\rt)^{C\ve V(\r)}.\eqno{(3.45)}$$
Substituting (3.43)-(3.45) into (3.42) yields
$$
|\wt\sq_g\p^n u|\le \ds\f{C\ve}{(1+t)(1+|\r|)^\nu}\ds\sum_{m=n+1}|\p^m u|+\ds\f{C\ve^2}{(1+t)^2(1+|\r|)^{n+2\nu}}\lt(\ds\f{1+t}{1+|\r|}\rt)^{C\ve V(\r)}.\eqno{(3.46)}
$$
On the other hand, by (2.6), (2.39) and (3.1), one has
$$|Z^I\p^n u|\le C(1+|q|)^{-n}\ds\sum_{J\le n+I}|Z^J u|\le C\ve(1+|q|)^{-n}(1+t)^{-\nu}\le C\ve (1+t)^{-\nu+nC\ve}(1+|\r|)^{-n-nC\ve}.\eqno{(3.47)}$$
In terms of (3.46) and (3.47), we know that both the assumptions in Lemma 6.6 of [20] hold. Thus, it follows
from Lemma 6.6  of [20] that (3.40) holds for $k=n+1$. And the proof of Lemma 3.7 is completed.
\qquad \qquad $\square$

\vskip 0.2 true cm
Finally, we derive the decay estimate of $|\p^k Z^I u|$ with $\max(1,k)+I\le N-4$.
\vskip 0.2 true cm
{\bf Lemma 3.8.} {\it Assume that the solution $u$ to (1.3) satisfies (3.1)-(3.3). For $\max(1,k)+I\le N-4$, then we have
$$|\p^k Z^I u|\le C(1+t)^{-1+C\ve}(1+|q|)^{1-k-\nu}.\eqno{(3.48)}$$}

{\bf Proof.} Thanks to (2.39), we just need to show
$$|\p^k Z^I u|\le C(1+t)^{-1+C\ve}(1+|\r|)^{1-k-\nu}.\eqno{(3.49)}$$
Next we use the induction method to prove (3.49).

For $I=0$ and all $k$, we have already proved (3.49) in (3.40) of Lemma 3.7.

Assume that (3.49) holds for $I\le m-1$ ($m\ge 1$)  and  all $k$, then we can prove (3.49) for $I=m$ and $k\le 1$.
In fact, by (3.20), we have for  $I=m$
\begin{align*}
|\wt\sq_g Z^m u|\le &C|Z^m u||\p^2 u|+C \ds\sum_{J+K\le m, J\le m-1, K\le m-1}|Z^J u||\p^2 Z^K u|\\
&+ \ds\sum_{J+K\le m, K\le m-1}|G(Z^J u,Z^K u)|.\tag{3.50}
\end{align*}
Hence, by (3.24) applied to $\phi=Z^m u$ and (3.50), together with (3.5), (3.27) and (3.49) for $m-1$, one has
\begin{align*}
|L^i_2&\p_i(r\p_\r Z^m u)|\le \ds\f{Cr}{\r_q}\bigg(|Z^m u||\p^2 u|+\ds\sum_{J+K\le m,J\le m-1,K\le m-1}|Z^J u||\p^2 Z^K u|\\
&\quad +\ds\sum_{J+K\le m,K\le m-1}|G(Z^J u,Z^K u)|\bigg)+C\ve|\p_\r Z^m u|+\ds\f{C\r^{-1}_q}{1+t}\ds\sum_{J\le m+2}|Z^J u|\\
&\le C\ve(\ds\f{|Z^m u|}{1+|\r|}+|\p_\r Z^m u|)+\ds\f{C\ve(1+|\r|)^{-C\ve}}{(1+t)^{1+\nu-C\ve}}\\
&\q+\ds\f{C\ve^2}{(1+t)^{1-C\ve}(1+|\rho|)^{2\nu}}+\ds\f{Cr}{\r_q}\ds\sum_{J+K\le m, K\le m-1}|G(Z^J u,Z^K u)|).\tag{3.51}
\end{align*}
Since
$$\ds\sum_{J+K\le m,K\le m-1}|G(Z^J u,Z^K u)|=|G(Z^m u,u)|+\ds\sum_{J\le m-1,K\le m-1}|G(Z^J u,Z^K u)|,$$
by Lemma 2.5, (3.1), (2.38), (3.27) or (3.49) for $m-1$, we have respectively
\begin{align*}
Cr\r^{-1}_q|G(Z^m u,u)|&\le Cr\r^{-1}_q(1+t)^{-1}(|Z^{m+1} u||\p^2 u|+|\p Z^m u||Z\p u|)\\
&\le \ds\f{C\ve^2}{(1+t)^{1-C\ve}(1+|\rho|)^{2\nu}}\tag{3.52}
\end{align*}
and
$$Cr\r^{-1}_q \ds\sum_{J\le m-1,K\le m-1}|G(Z^J u,Z^K u)|\le
\ds\f{C\ve^2}{(1+t)^{1-C\ve}(1+|\rho|)^{2\nu}}.\eqno{(3.53)}$$
Substituting (3.52)-(3.53) into (3.51) yields
$$
|L^i_2\p_i(r\p_\r Z^m u)|
\le C\ve\biggl(\ds\f{|Z^m u|}{1+|\r|}+|\p_\r Z^m u|\biggr)+\ds\f{C\ve(1+|\r|)^{-C\ve}}{(1+t)^{1+\nu-C\ve}}
+\ds\f{C\ve^2}{(1+t)^{1-C\ve}(1+|\rho|)^{2\nu}}.\eqno{(3.54)}
$$
Then it follows from (3.54) and Lemma 3.5 that
$$|Z^m u|(1+|\rho|)^{-1}+|\p_\rho Z^m u|\le C\ve (1+t)^{-1+C\ve}(1+|\rho|)^{-\nu},\eqno{(3.55)}$$
which means that (3.49) holds for $I=m$ and $k\le 1$.

Finally, assuming (3.49) for  $I\le m$ and all $k\le n$, and (3.49) for $I\le m-1$ and all $k$,
we now show that (3.49) holds for $I=m\ge 1$ and $k=n+1\ge 2$.

In this case, by (3.34)-(3.35) together with (3.49) for $k\le n$ and $I\le m $, and $k\le n+2$ and $I\le m-1$, one has
\begin{align*}
|\wt\sq_g\p^n Z^I u|&\le C|\p u|\ds\sum_{m=n+1}|\p^m Z^I u|+\ds\sum_{J+K\le I, k\le n, m+k=n+2,
\text{$m\le n$ or $k\le I-1$}}|\p^m Z^J u||\p^k Z^K u |\\
&\;\;\;\;+\sum_{m+k=n, J+K\le I, k+K<n+I}|G(\p^m Z^J u,\p^k Z^K u)|\\
&\le \ds\f{C\ve}{(1+t)(1+|\r|)^\nu}\ds\sum_{m=n+1}|\p^m Z^I u|+\ds\f{C\ve^2}{(1+t)^{2-C\ve}(1+|\rho|)^{n+2\nu}},\tag{3.56}
\end{align*}
here we point out that $\ds\sum_{m+k=n, J+K\le I, k+K<n+I}|G(\p^m Z^J u,\p^k Z^K u)|$ has been treated as in (3.43)-(3.45).
From (3.56), as shown in Lemma 3.7, we can show that (3.49) holds for $I=m\ge 1$ and $k=n+1\ge 2$. By induction
method and all the analysis above, we complete the proof of Lemma 3.8.\qquad $\square$

In summary, collecting Lemma 3.1, Lemma 3.4 and Lemma 3.6-Lemma 3.8, we arrive at

{\bf Proposition 3.9.} {\it Assume that $u$ is the solution to (1.3) and (3.1)-(3.3) hold.
Then
\begin{equation*}
|u|\le\left\{
\begin{aligned}
& C \ve (1+t)^{-1}(1+|q|), \\
& C \ve (1+t)^{-1+C \ve}(1+|\rho|)^{1-\nu-C\ve},
\end{aligned}
\right.\tag{3.57}
\end{equation*}
$$|\p u|\le C\ve (1+t)^{-1}(1+|\r|)^{-\nu},\eqno{(3.58)}$$
$$|\p^2 u|\le C\ve (1+t)^{-1}(1+|\r|)^{-1-\nu}\r_q,\eqno{(3.59)}$$
\begin{equation*}
|Zu|\le\left\{
\begin{aligned}
&C \ve (1+t)^{-1+C\ve}(1+|\rho|)^{1-\nu}, \\
&C\ve (1+t)^{-1}(|q|+(1+t)^{C\ve}),
\end{aligned}
\right.\tag{3.60}
\end{equation*}
furthermore,
$$|\p^k Z^I u|\le C\ve (1+t)^{-1+C\ve}(1+|q|)^{1-k-\nu}\qquad\text{for $\max(1, k)+I\le N-4$.}\eqno{(3.61)}$$}

\vskip 0.4 true cm
\centerline {\bf $\S 4$. Weighted energy estimates for the nonlinear problem (1.3)}
\vskip 0.4 true cm

As in [4] and [20], we now establish the weighted energy estimates for
the equation (1.3) so that the standard energy $E_N(t)=\ds\sum_{I\le N}\int |\p Z^I u(t,x)|^2 dx\le C\ve^2(1+t)^{C\ve}$
can be shown under the weak assumption $E_N(t)\le  C\ve^2(1+t)^{\dl}$ for some small fixed constant $0<\dl<\f12$.
From this and the higher order energy estimates in $\S 5$, the validity of the weak decay estimate (3.1) of $u$ can be proved in $\S 6$.

We  will  choose such a weight in the weighted energy of (1.3)
$$W=e^{\si(t) V(\r)+\vp(q)}\eqno{(4.1)}$$
with
\begin{equation*}
\left\{
\begin{aligned}
&\si(t)=\k\ve\ln|1+t|,\\
& V(\rho)={|\rho-2|}^{-\nu'}, \\
&\vp'(q)=(1+|q|)^{-3/2},\\
\end{aligned}
\right.
\end{equation*}
where $\k>0$ and $\nu'>1/2$ are constants, the function $\r(t,x)$ has been defined in (2.34),
and $q=r-t$. In addition, $\r\le 1$ is known.

Now we state the weighted energy estimate in this section.
\vskip 0.2 true cm
{\bf Proposition 4.1.} {\it Assume that $u$ is a solution to (1.3) and the notations $g, c, E$ (corresponding to $g^{ij}, c^{ij}, E^{ij}$,
respectively) have
been defined in previous sections.  In addition, the following assumptions are fulfilled
\begin{align*}
&|g-c-E|\le \ds\f{1}{4},\qquad |E|<\ds\f{1}{8}, \qquad |\p(g-E)|\le \ds\f{C\ve}{1+t},\tag{4.2}\\
&(1+|q|)^{-1}|Z\p u|+|\p^2 u|\le C\ve(1+t)^{-1+C\ve}(1+|q|)^{-1-\nu'},\tag{4.3}\\
&\r_t<0, \quad \ds\f{g^{ij}\r_i\r_j}{\r_t(1+|\r|)^{1+\nu'}}\ge -\ds\f{1}{\k\ve(1+t)\ln(1+t)},\tag{4.4}\\
\end{align*}
where $\k$ and $\nu'>1/2$ are given in (4.1),
then we have
\begin{align*}
&\ds\int_{\Sigma_t}|\p u|^2 W dx +\int^t_0\int_{\Sigma_\tau}(1+|q|)^{-3/2}|\bar Z u|^2 W dxd\tau\\
&\le C\int_{\Sigma_0}|\p u|^2 W dx
+\int^t_0\int_{\Sigma_\tau}\ds\f{C\ve}{1+\tau}|\p u|^2 W dxd\tau +\int^t_0\int_{\Sigma_\tau}\ds\f{C(1+\tau)}{\ve}|\wt\sq_g u|^2 W dxd\tau,\tag{4.5}
\end{align*}
where $W$ is defined in (4.1), $\Sigma_t=\{x\in\Bbb R^3: |x|\le 1+t\}$,
and $|\bar Z u|^2=\ds\sum^3_{\al=1}(\p_\al u+\o_\al\p_t u)^2$.
}

\vskip 0.2 true cm

{\bf Remark 4.1.} {\it The choice of the complicated weight function $W$ in (4.1) is due to the following two reasons: First, the factor $e^{\si(t) V(\r)}$
in $W$ comes from the weight in [4] and [20] which is used to derive the weighted energy estimate for the wave equation $\ds\sum_{i,j=0}^3g^{ij}(u)\p_{ij}^2u=0$.
Motivated by this, we intend to search a new weight $W$ containing the factor $e^{\si(t) V(\r)}$
so that  the  weighted energy estimate for the wave
equation $\ds\sum_{i,j=0}^3g^{ij}(u, \p u)\p_{ij}^2u=0$ can be derived. Second, due to the appearance of the term $e_{k}^{ij}\p_ku\p_{ij}^2u$
in (1.3) and the related null condition property, we want to add another factor $e^{\vp(q)}$ in the new weight $W$.
In this case, we may obtain the good controls for the troublesome terms $Q_2\p_t^2u$ and $I_1$ in (4.17) and (4.19) below,
which are induced by the appearances of $\ds\sum_{k=0}^3 e^{ij}_k\p_ku$ in the coefficients $g^{ij}(u,\p u)$. The so-called ``good control"
means that the corresponding integral can be finally absorbed by the left hand side of (4.5).}

\vskip 0.2 true cm
{\bf Remark 4.2.} {\it By the way, if the null condition does not
hold (for this case, we have actually proved the blowup of solution to (1.3) in [7-8]), that is, $\ds\sum^3_{k=0}\ds\f{1}{2}e^{00}_k\o_k-e^{\b 0}_k\o_\b\o_k+\ds\f{1}{2}e^{\al\b}_k\o_\al\o_\b\o_k \not\equiv 0$ holds, then we can not obtain the energy estimate (4.5)
since the troublesome term $ |u^2_t(\ds\f{1}{2}e^{00}_k\o_k-e^{\b 0}_k\o_\b\o_k+\ds\f{1}{2}e^{\al\b}_k\o_\al\o_\b\o_k) \p^2_t u|\lt(\le C\ve(1+\tau)^{-1+C\ve} |\p u|^2 \text{ by the assumption }(4.3) \rt)$
will appear in the right hand side of $|I_1|$ in (4.19) below,
 which is not a good control
term (only $C\ve(1+\tau)^{-1} |\p u|^2 $ is a good control term by Gronwall's inequality).}

\vskip 0.2 true cm

{\bf Proof.} Denote by $v_{\al}=\p_{\al} v$ for $\al =1,2,3$ and $v_t=\p_t v=v_0$, then it follows from a direct computation that
\begin{align*}
W\p_t u\wt\sq_g u&=(g^{00}\p^2_t u\p_t u+2g^{0\al}\p_0\p_\al u\p_t u+g^{\al\b}\p_\al\p_\b u\p_t u)W\\
&=\ds\f{1}{2}\p_t\lt[(g^{00}u^2_t-g^{\al\b}u_\al u_\b)W\rt]-\ds\f{1}{2}(g^{00}u^2_t-g^{\al\b}u_\al u_\b)\p_t W\\
&\;\;\;\;+\p_\al\lt[(g^{0\al}u^2_t+g^{\al\b}u_t u_\b)W\rt]-(g^{0\al}u^2_t+g^{\al\b}u_t u_\b)\p_\al W\\
&\;\;\;\;+W\lt(-\ds\f{1}{2}u^2_t\p_t g^{00}-u^2_t\p_\al g^{0\al}-u_\b u_t\p_\al g^{\al\b}+\ds\f{1}{2}u_\al u_\b\p_t g^{\al\b}\rt).\tag{4.6}
\end{align*}
Integrating (4.6) over $x\in\Bbb R^3$ yields
$$\ds\f{1}{2}\ds\f{d}{dt}\int(-g^{00}u^2_t+g^{\al\b}u_\al u_\b)W dx =-\int u_t W\wt\sq_g u dx+\int I_1 W dx+\int I_2 dx,\eqno{(4.7)}$$
where
\begin{align*}
I_1&=-\ds\f{1}{2}u^2_t\p_t g^{00}-u^2_t\p_\al g^{0\al}-u_\b u_t\p_\al g^{\al\b}+\ds\f{1}{2}u_\al u_\b\p_t g^{\al\b}\\
&=\ds\f{1}{2}u_i u_j\p_t g^{ij}-u_t u_j\p_i g^{ij},\\
I_2&=-\ds\f{1}{2}(g^{00}u^2_t-g^{\al\b}u_\al u_\b)\p_t W-(g^{0\al}u^2_t+g^{\al\b}u_t u_\b)\p_\al W\\
&=\ds\f{1}{2}g^{ij}u_i u_j W_t-g^{ij}u_i u_t W_j.
\end{align*}
Next we analyze each term in the right hand side of (4.7).
At first, we treat $I_1$ in (4.7). Since $g^{ij}=c^{ij}+D^{ij}+E^{ij}$  and  $|\p(g-E)|\le \ds\f{C\ve}{1+t}$, it follows that
$$ |\ds\f{1}{2}u_i u_j\p_t (c^{ij}+D^{ij})-u_t u_j\p_i (c^{ij}+D^{ij} ) |\le C\ve (1+t)^{-1}|\p u|^2.\eqno{(4.8)}$$
In addition, by a direct computation, we have
\begin{align*}
&\ds\f{1}{2}u_i u_j\p_t E^{ij}-u_t u_j\p_i E^{ij}\\
=\;\;&\ds\f{1}{2}u_i u_j\ds\sum^3_{k=0}e^{ij}_k\p_t\p_k u-u_t u_j\ds\sum ^3_{k=0}e^{ij}_k\p_i\p_k u\\
=\;\;&u^2_t \ds\sum ^3_{k=0} (-\ds\f{1}{2}e^{00}_k\p_t\p_k u-e^{\b 0}_k\p_\b\p_k u)+u_t u_\b(-\ds\sum ^3_{k=0}e^{\al\b}_k\p_\al\p_k u)
+\ds\f{1}{2}u_\al u_\b \ds\sum ^3_{k=0} e^{\al\b}_k\p_t\p_k u.\tag{4.9}
\end{align*}
Substituting $\p_t\p_k u=\p_t\p_k u+\o_k\p^2_t u-\o_k\p^2_t u$  and  $\p_\b\p_k u=\p_\b\p_k u-\o_\b\o_k\p^2_t u+\o_\b\o_k\p^2_t u$
into the first term in the right hand side of (4.9) yields
\begin{align*}
&u^2_t \ds\sum ^3_{k=0}(-\ds\f{1}{2}e^{00}_k\p_t\p_k u-e^{\b 0}_k\p_\b\p_k u)\\
&=u^2_t \ds\sum ^3_{k=0} \lt[-\ds\f{1}{2}e^{00}_k(\p_t\p_k u+\o_k\p^2_t u)-e^{\b 0}_k (\p_\b\p_k u-\o_\b\o_k\p^2_t u)\rt]+u^2_t\ds\sum^3_{k=0}(\ds\f{1}{2}e^{00}_k\o_k-e^{\b 0}_k\o_\b\o_k)\p^2_t u\\
&=-e^{\b0}_0 u^2_t(\p_\b\p_tu+\o_\b\p^2_tu)+u^2_t \ds\sum ^3_{k=1} \lt[-\ds\f{1}{2}e^{00}_k(\p_t\p_k u+\o_k\p^2_t u)-e^{\b 0}_k (\p_\b\p_k u-\o_\b\o_k\p^2_t u)\rt]\\
&\q\q+u^2_t\ds\sum^3_{k=0}(\ds\f{1}{2}e^{00}_k\o_k-e^{\b 0}_k\o_\b\o_k)\p^2_t u.\tag{4.10}
\end{align*}
Similarly, the last two terms in the right hand side of (4.9) admit the following expressions respectively
\begin{align*}
&u_t u_\b(-\ds\sum ^3_{k=0}e^{\al\b}_k\p_\al\p_k u)=-u_t u_\b e^{\al\b}_0(\p_\al\p_t u+\o_\al\p^2_t u)-u_t u_\b \ds\sum ^3_{k=1}e^{\al\b}_k(\p_\al\p_k u-\o_\al\o_k\p^2_t u)\\
&\qquad \qquad \qquad \qquad \qquad \qquad \qquad -u_tu_\b \ds\sum ^3_{k=0}e^{\al\b}_k\o_\al\o_k\p_t^2u,\tag{4.11}\\
&\ds\f{1}{2}u_\al u_\b \ds\sum ^3_{k=0}e^{\al\b}_k\p_t\p_k u=\ds\f{1}{2}u_\al u_\b \ds\sum ^3_{k=1}e^{\al\b}_k(\p_t\p_k u+\o_k\p^2_t u)-\ds\f{1}{2}u_\al u_\b \ds\sum ^3_{k=0}e^{\al\b}_k
\o_k\p^2_t u.\tag{4.12}
\end{align*}
Note that for $k=1,2,3$ and by (4.3),
$$|\p_t\p_k u+\o_k\p^2_t u|=|t^{-1}\G_{0k}\p_t u-\o_k t^{-1}q\p^2_t u|\le C(1+t)^{-1}|Z\p u|\le C\ve(1+t)^{-1}\eqno{(4.13)}$$
and
$$|\p_\b \p_k u-\o_\b \o_k\p^2_t u|=|t^{-1}(\G_{0 \b}\p_k u -\o_\b\G_{0 k}\p_t u)-t^{-1}\o_\b q(\p_t\p_k u-\o_k\p^2_t u)|\le C\ve(1+t)^{-1}.\eqno{(4.14)}$$
Then substituting (4.10)-(4.14) into (4.9) derives
$$\ds\f{1}{2}u_i u_j\p_t E^{ij}-u_t u_j\p_i E^{ij}=Q_1+Q_2\p^2_t u,\eqno{(4.15)}$$
where $ |Q_1|\le C\ve(1+t)^{-1}|\p u|^2$ and $Q_2=\ds\sum^3_{k=0}u^2_t(\ds\f{1}{2}e^{00}_k\o_k-e^{\b 0}_k\o_\b\o_k)-u_t u_\b e^{\al\b}_k\o_\al\o_k-\ds\f{1}{2}u_\al u_\b e^{\al\b}_k \o_k$. Here we point out that the term $Q_2\p^2_t u$ should be specially treated as follows.

Set
$ X_\al=u_\al+\o_\al u_t$, then $u_\al=X_\al -\o_\al u_t$.
Thanks to  the null condition property, we can obtain
\begin{align*}
Q_2&=\ds\sum^3_{k=0}u^2_t(\ds\f{1}{2}e^{00}_k\o_k-e^{\b 0}_k\o_\b\o_k)-u_t (X_\b -\o_\b u_t) e^{\al\b}_k\o_\al\o_k-\ds\f{1}{2}(X_\al -\o_\al u_t) (X_\b -\o_\b u_t) e^{\al\b}_k \o_k\\
&=\ds\sum^3_{k=0}u^2_t(\ds\f{1}{2}e^{00}_k\o_k-e^{\b 0}_k\o_\b\o_k+\ds\f{1}{2}e^{\al\b}_k\o_\al\o_\b\o_k)-\ds\f{1}{2}X_\al X_\b e^{\al\b}_k\o_k\\
&=-\ds\f{1}{2}\ds\sum^3_{k=0}X_\al X_\b e^{\al\b}_k\o_k\qquad \text{(due to $\ds\sum^3_{k=0}\ds\f{1}{2}e^{00}_k\o_k-e^{\b 0}_k\o_\b\o_k+\ds\f{1}{2}e^{\al\b}_k\o_\al\o_\b\o_k\equiv 0$)}.\tag{4.16}
\end{align*}
Thus, by (4.16) and (4.3), one has
$$
|Q_2 \p^2_t u|=\ds\f12|\ds\sum^3_{k=0}X_\al X_\b e^{\al\b}_k\o_k \p^2_t u|\le C \ds\sum^3_{\al=1}X^2_\al |\p^2_t u|
\le C\ve(1+t)^{-1+C\ve}(1+|q|)^{-1-\nu'}|\bar Z u|^2.\eqno{(4.17)}
$$
Obviously, if the null condition does not hold here, then $\ds\sum^3_{k=0}\ds\f{1}{2}e^{00}_k\o_k-e^{\b 0}_k\o_\b\o_k+\ds\f{1}{2}e^{\al\b}_k\o_\al\o_\b\o_k
\not\equiv 0$ and we only obtain
$$|Q_2 \p^2_t u|\le  C\ve(1+t)^{-1+C\ve}(1+|q|)^{-1-\nu'}|\p u|^2+C\ve(1+t)^{-1+C\ve}(1+|q|)^{-1-\nu'}|\bar Z u|^2,\eqno{(4.18)}$$
which means that the first term in the right hand side of (4.18) can not be absorbed globally by the term $\ds\int_{\Sigma_t}|\p u|^2 W dx$
of (4.5) (the concrete reason is: in this case, (4.5) becomes $\int_{\Sigma_t}|\p u|^2Wdx+...\le \int_0^t\int_{\Sigma_\tau}
C\ve(1+t)^{-1+C\ve}|\p u|^2Wdxd\tau+...$, then the crucial integral $\int_{\Sigma_t}|\p u|^2Wdx$ cannot be uniformly
controlled by Gronwall's inequality since $\int_0^{\infty}(1+t)^{-1+C\ve}dt=\infty$).

Consequently, collecting (4.8), (4.15) and (4.17) yields
$$|I_1|\le C\ve(1+t)^{-1}|\p u|^2+C\ve(1+t)^{-1+C\ve}(1+|q|)^{-1-\nu'}|\bar Z u|^2.\eqno{(4.19)}$$

Next let us deal with the term $I_2$ in (4.7). Set
$$W=\wt W e^{\vp(q)},  \q\text {where} \;\;\wt W=e^{\si(t)V(\r)}.$$
Then it follows from a direct computation that
\begin{align*}
I_2&=(\ds\f{1}{2}g^{ij}u_i u_j-g^{i0}u_i u_t )\bigg((A\r_t +B)W-\vp'(q)W\bigg)-g^{i\al}u_i u_t\bigg(A\r_\al W+\vp'(q)\o_\al W \bigg)\\
&=W(I_{21}+I_{22}+I_{23})\tag{4.20}
\end{align*}
with
\begin{align*}
&I_{21}=\lt((\ds\f{1}{2}g^{ij}u_i u_j-g^{i0}u_i u_t)\r_t-g^{i\al}u_i u_t\r_\al\rt)A,\\
&I_{22}=\lt(\ds\f{1}{2}g^{ij}u_i u_j-g^{i0}u_i u_t\rt)B,\\
&I_{23}=\vp'(q)\lt(-\ds\f{1}{2}g^{ij}u_i u_j+g^{i0}u_i u_t-g^{i\al}u_i u_t\o_\al\rt),
\end{align*}
where $A=\ds\f{\k\nu'\ve\ln|1+t|}{|\r-2|^{1+\nu'}}$ and $B=\ds\f{\k\ve}{(1+t)|\r-2|^{\nu'}}$. We now treat
each term in the right hand side of (4.20).

Since $|g-c-E|<1/4$ and  $|E|<1/8$,  we have $|g-c|<1/2$,  which means that
the $3\times 3$ matrix $({g^{\al\b}})_{\al,\beta=1}^3$ is nonnegative definite.
Thus, as in (7.9) of [20], we have by (4.4)

$$I_{21}\le -\ds\f{1}{2}A\ds\f{g^{ij}\r_i\r_j u^2_t}{\r_t}
\le\f{C\ve}{1+t}|\p u|^2.\eqno{(4.21)}$$
On the other hand, a direct computation yields
$$\ds\f{1}{2}g^{ij}u_i u_j-g^{i0}u_i u_t=\ds\f{1}{2}(-g^{00}u^2_t+g^{\al\b}u_\al u_\b).\eqno{(4.22)}$$
In addition,  it follows from  $|D|<1/4$ and $|E|<1/8$   that
$$\ds\f{1}{2}|\p u|^2=\ds\f{1}{2}(u^2_t +\dl^{\al\b}u_\al u_\b)\le -g^{00}u^2_t
+g^{\al\b}u_\al u_\b \le 2(u^2_t +\dl^{\al\b}u_\al u_\b)=2|\p u|^2.\eqno{(4.23)}$$
This, together with (4.22) and $\r\le 1$, yields
$$I_{22}\le \ds\f{C\ve}{1+t}|\p u|^2.\eqno{(4.24)}$$
We now treat the third term in (4.20). It is noted that
\begin{align*}
-\ds\f{1}{2}c^{ij}u_i u_j+c^{i0}u_i u_t-c^{i\al}u_i u_t\o_\al&=-\ds\f{1}{2}(-u^2_t +\ds\sum^3_{\al=1}u^2_\al)-u^2_t-u_\al u_t\o_\al\\
&=-\ds\f{1}{2}\ds\sum^3_{\al=1}(u_\al+\o_\al u_t)^2=-\ds\f{1}{2}|\bar Z u|^2.\tag{4.25}
\end{align*}
Additionally, due to  $\vp '(q)=(1+|q|)^{-3/2}$, it follows that
\begin{align*}
&\lt|\vp'(q)\lt[-\ds\f{1}{2}(D^{ij}+E^{ij})u_i u_j+(D^{i0}+E^{i0})u_i u_t-(D^{i\al}+E^{i\al})u_i u_t\o_\al\rt]\rt|\\
\le \;\;&C|\vp'(q)|(|u|+|\p u|)|\p u|^2\\
\le \;\;&C\ve(1+t)^{-1}(1+|q|)|\vp'(q)||\p u|^2\\
\le \;\;&C\ve(1+t)^{-1}|\p u|^2.
\end{align*}
This, together with (4.25), yields
$$I_{23}\le -\ds\f12(1+|q|)^{-\f32}|\bar Z u|^2+C\ve(1+t)^{-1}|\p u|^2.\eqno{(4.26)}$$
Thus, by substituting (4.21), (4.24) and (4.26) into (4.20), we have
$$I_2\le \biggl(-\ds\f{1}{2}(1+|q|)^{-\f32}|\bar Zu|^2+C\ve(1+t)^{-1}|\p u|^2\biggr)W.\eqno{(4.27)}$$
Noting that
\begin{align*}
\int^t_0\int(-g^{ij}\p_i\p_j u)u_t W dxd\tau&=\int^t_0\int -u_t\wt\sq_g u Wdxd\tau\\
&\le \int^t_0\ds\f{C\ve}{1+\tau}\int_{\Sigma_\tau}|\p u|^2 W dxd\tau+\ds\f{C}{\ve}\int^t_0\int_{\Sigma_\tau}(1+\tau)|\wt\sq_g u|^2 Wdxd\tau,\tag{4.28}
\end{align*}
then integrating (4.7) over $[0, t]$ and applying (4.19), (4.23) and (4.27)-(4.28), we can arrive at
\begin{align*}
&\ds\f{1}{4}\int_{\Sigma_t}|\p u|^2 Wdx+\ds\f12\int^t_0\int_{\Sigma_\tau}\lt[1-C\ve(1+t)^{-1+C\ve}(1+|q|)^{-\nu'+1/2}\rt](1+|q|)^{-3/2}|\bar Z u|^2 Wdxd\tau\\
\le\;& \ds\f12\int_{\Sigma_0}(-g^{00}u^2_t+g^{\al\b}u_\al u_\b)W dx+\ds\f{C}{\ve}\int^t_0\int_{\Sigma_\tau}(1+\tau)|\wt\sq_g u|^2 Wdxd\tau
+\int^t_0\int_{\Sigma_\tau}\ds\f{C\ve}{1+\tau}|\p u|^2 Wdx\\
\le\;& \int_{\Sigma_0}|\p u|^2 Wdx+\ds\f{C}{\ve}\int^t_0\int_{\Sigma_\tau}(1+\tau)|\wt\sq_g u|^2 Wdxd\tau
+\int^t_0\int_{\Sigma_\tau}\ds\f{C\ve}{1+\tau}|\p u|^2 Wdx.\tag{4.29}
\end{align*}
Due to $\nu'>\f12$  and $1-C\ve (1+t)^{-1+C\ve} \ge 1/2$ for small $\ve>0$,  we then obtain (4.5)
from (4.29).\qquad $\square$

\vskip 0.4 true cm
\centerline {\bf $\S 5$. Higher order energy estimates for the problem (1.3)}
\vskip 0.4 true cm

In this section, under the strong decay  assumptions on the solution $u$ to (1.3), which are given in Prop.3.9,
we will give the higher order energy estimates for $E_{k,i}(t)=\ds\sum_{0\le \text{\bf k}\le k , 0\le I\le i} \int |\p \p^{\text{\bf k}} Z^I u|^2 W dx$. Before doing this,
we  show a weighted Poincar\'e lemma similar to that in Lemma 8.1 of [20].

{\bf Lemma 5.1.} {\it Assume that $W$ is defined in (4.1) with suitably large $\k>0$ and
$$|\p_{\r}\p_r \r|\le \ds\f{C\ve\ln|1+t|}{(1+|\r|)^{1+\nu'}}\p_r\r,\;\;\;0<\p_r\r<\infty.\eqno{(5.1)}$$
Then for function $u$ supported in $r\le 1+t$ ,
$$\int \lt(\ds\f{|u|}{1+|\r|}\p_r\r\rt)^2 W dx+\int \lt(\ds\f{|u|}{1+|r-t|}\rt)^2 W dx\le C\int |\p u|^2 W dx.\eqno{(5.2)}$$}

{\bf Proof.} Although the proof of (5.1) is completely similar to that in Lemma 8.1 of [20], we will still give the
details for reader's convenience, due to the different form of the weight $W$.
Notice that we only require to treat the first integral in the left hand side of (5.2) since the second one is a special case of the first with $\rho=r-t$.

Taking integration by parts, we have
\begin{align*}
&\int ^\infty_0 \lt(\ds\f{|u|}{|\rho -2|}\p_r\r\rt)^2 W r^2 dr= \int^1_{-\infty}\lt(\ds\f{|u|}{|\rho -2|}\rt)^2\p_r\r W r^2 d\rho
= \int^1_{-\infty} |u|^2\p_r\r W r^2 d\lt(\ds\f{1}{|\rho -2|}\rt)\\
&\qquad =-2\int^1_{-\infty} \ds\f{u}{|\rho -2|}\p_{\rho} u\p_r\r W r^2 d\rho-\int^1_{-\infty}\lt(\ds\f{u}{|\rho-2|}\rt)^2|\rho-2|
\p_\rho\lt(\p_r\r W r^2\rt)d\rho.\tag{5.3}
\end{align*}
In addition, by (5.1), $1+|\rho|\ge \ds\f12 |\rho-2|$ and  suitably large $\k$, we have
\begin{align*}
\p_\rho&\lt(\p_r\r  W r^2\rt)=\lt(\p_{\rho} \p_r\r \rt)W r^2+\p_r\r \lt[(\p_\rho W) r^2+2rW \p_{\rho} r \rt]\\
&\ge -\ds\f{C\ve\ln |1+t|}{(1+|\rho|)^{1+\nu'}}(\p_r \rho)W r^2 +\p_r\r \left[\ds\f{\k \nu'\ve\ln |1+t|}{|\rho-2|^{1+\nu'}}+(1+|q|)^{-3/2}(\rho_q)^{-1}\right]W r^2+2rW\\
&\ge(1+|q|)^{-3/2}(\rho_q)^{-1} r^2 W \p_r {\rho}+2rW\ge 0.\\
\end{align*}
This, together with (5.3) and  H\"older inequality, yields
$$\int ^\infty _0 \lt(\ds\f{|u|}{|\rho -2|}\p_r\r \rt)^2W r^2 dr\le 2\lt(\int^\infty_0 \lt(\ds\f{|u|}{|\rho -2|}\p_r\r \rt)^2 W r^2 dr\rt)^{1/2}\lt(\int^\infty_0 \lt(\p_{\rho} u\p_r\r \rt)^2 W r^2 dr\rt)^{1/2}$$
and further
$$\int ^\infty _0 \lt(\ds\f{|u|}{|\rho -2|}\p_r\r \rt)^2W r^2 dr\le 4\int^\infty_0 \lt(\p_{\rho} u\p_r\r \rt)^2 W r^2 dr=4\int^\infty_0 \lt(\p_r u\rt)^2W r^2 dr\eqno{(5.4)}$$
Integrating (5.4) over the angular variables and using $\ds\f12|\rho-2|\le 1+|\rho|$, we obtain (5.2).\qquad$\square$

For the energy
$$ E_{k,i}(t)=\ds\sum_{0\le I\le i,0\le\text{\bf k}\le k} \int |\p \p^{\text{\bf k}} Z^I u|^2 W dx\eqno{(5.5)}$$
where $W$ is defined in (4.1) with suitably large $\k$,  we have

{\bf Proposition 5.2.} {\it Let $N\ge 8$ be fixed, $1/2<\nu'<1$ and $N'=[N/2]+2$. Assume that $u$ is a solution of (1.3) for $0\le t< T$ and
$$|\p u|\le \ds\f{C\ve}{1+t}\ds\f{1}{(1+|\r|)^{\nu'}},\eqno{(5.6)}$$
$$|\p^2 u|\le \ds\f{C\ve\r_q}{1+t}\ds\f{1}{(1+|\r|)^{1+\nu'}},\eqno{(5.7)}$$
$$|u|+|Z u|\le C\ve (1+t)^{-1+C\ve}(1+|\r|)^{1-\nu'},\eqno{(5.8)}$$
$$|\p Z^I u|+(1+|q|)^{-1}|Z^I u|\le\ds\f{C\ve}{1+t}(1+t)^{C\ve}(1+|q|)^{-\nu'}\;\;\;for \quad I\le N',\eqno{(5.9)}$$
Then  for $k+i\le N,\;\; k, i\ge 0$,
$$E_{k,i}(t)\le C E_{k,i}(0)+\int^t_0 \ds\f{C\ve}{1+\tau}E_{k,i}(\tau) d\tau +\int^t_0\ds\f{C\ve}{(1+\tau)^{1-C \ve}}(E_{k+1,i-1}(\tau)+E_{k-1,i}(\tau))d\tau, \eqno{(5.10)}$$
where $E_{-1,n}=0$ and $E_{m,-1}=0$.}
\vskip 0.3 true cm

{\bf Proof.}  We will use Proposition 4.1 and Lemma 5.1 to prove (5.10). To this end, we need to  verify  all the  assumptions (4.2)-(4.4)
of Proposition 4.1 and (5.1) of Lemma 5.1 respectively.

First, let us verify the assumptions (4.2) of Proposition 4.1. By applying (5.8) and  the facts of $ |\r|\le Ct$  and  $\nu' >1/2$, we have
$$|u|+|Z u|\le C\ve (1+t)^{-1+C\ve}(1+|\r|)^{1-\nu'}\le C\ve (1+t)^{-1/2+C\ve}.$$
From this and (5.6), it is enough to assume $|D|\le \ds\f{1}{4}$  and  $|E|\le \ds\f{1}{8}$ for
sufficiently small $\ve>0$. In addition, (5.6) derives $|\p(g-E)|\le\ds\f{C\ve}{1+t}$ directly.
Thus, (4.2) holds.

Second, it is obvious that (5.6) implies (2.37) holds. By Remark 2.3, we know that (2.39) is true. As in (3.11),  one has $|\p_q E_{LL}|
\le \ds\f{C\ve}{1+t}\ds\f{1}{(1+|\r|)^{\nu'}}$ by (5.9) and (2.39). This, together with (5.6),  yields (2.38) of Lemma 2.4.
Then it follows from this,  (5.7) and (5.9) that
$$(1+|q|)^{-1}|Z\p u|+|\p^2 u|\le C\ve(1+t)^{-1+C\ve}(1+|q|)^{-1-\nu'},$$
which means that (4.3) holds.

Third, by $L^i_2\p_i \r=0$, we have $\p_p \r=-\ds\f{1}{4}(D_{LL}+E_{LL})\p_q \r$ and
$$-\r_t=-(\r_p-\r_q)=(1+\ds\f{1}{4}D_{LL}+\ds\f{1}{4}E_{LL})\p_q \r >\ds\f{1}{2}\r_q\ge \ds\f{1}{2}(\ds\f{1+t}{1+|\r|})^{-C\ve}.$$

From (5.6) and (5.8)-(5.9), we know that the assumption in Lemma 5.3 of [20] hold.
Then as in (7.3) of [20], we can arrive at
$$\ds\f{g^{ij}\r_i\r_j}{\r_t(1+|\r|)^{1+\nu'}}\ge -\ds\f{1/{(\k\nu')}}{(1+t)\ln(1+t)},$$
which means that (4.4) holds.

Fourth, we verify (5.1) of Lemma 5.1. To this end, we intend to  use Lemma 2.4  to  derive  (5.1).
We now  verify the assumption (2.40) of Lemma 2.4.  Due to (5.7),
we only need to verify that $\p^2_q E_{LL}$ satisfies (2.40).
As in (3.10), by the null condition, we have
$$\p_q E_{LL}=-\ds\f{1}{2}{\un L}^m e^{ij}_0(\p_t\p_m u+\o_m\p^2_t u)L_i L_j
-\ds\f{1}{2}{\un L}^m\ds\sum_{\al=1}^3 e^{ij}_{\al}(\p_{\al}\p_m u-\o_{\al}\o_m \p^2_t u)L_i L_j.\eqno{(5.11)}$$
For the first term in the right hand side of (5.11),
\begin{align*}
&\lt|\p_q(-\ds\f{1}{2}{\un L}^m e^{ij}_0(\p_t\p_m u+\o_m\p^2_t u)L_i L_j)\rt|=\lt|-\ds\f{1}{2}{\un L}^m  e^{ij}_0L_i L_j\p_q  (\p_t\p_m u+\o_m\p^2_t u) \rt|\\
&=\lt|-\ds\f{1}{4}{\un L}^{\al} e^{ij}_0 L_i L_j (\p_r-\p_t)(t^{-1}\Gamma_{0\al}\p_t u-\o_\al t^{-1}q \p^2_t u)\rt|\\
&\le C t^{-1}(|Z\p^2 u|+|\p^2 u|)+Ct^{-2}|Z\p u|\\
&\le C t^{-1}(1+|q|)^{-1}\ds\sum_{0\le I\le 2}|\p Z^I u|+Ct^{-2}|Z\p u|.
\end{align*}
This, together with (5.9), (2.38)-(2.39) and the fact of $|q|\le Ct$, yields
\begin{align*}
\lt|\p_q(-\ds\f{1}{2}{\un L}^m e^{ij}_0(\p_t\p_m u+\o_m\p^2_t u)L_i L_j)\rt|
&\le C\ve (1+t)^{-1}(1+|q|)^{-1-\nu'}(1+t)^{-1+C\ve}\\
&\le C\ve (1+t)^{-1}(1+|\rho|)^{-1-\nu'} \r_q.\tag{5.12}
\end{align*}
The second term  in the right hand side of (5.11) proceed similarly. Hence,
we have from (5.11)-(5.12)
$$|\p^2_q E_{LL}|\le C\ve (1+t)^{-1}(1+|\rho|)^{-1-\nu'} \r_q,$$
which means that the assumption (2.40) in Lemma 2.4 holds. Then by (2.41) of Lemma 2.4, one has
$$|\p_{\r}\p_q \r|\le \ds\f{C\ve}{(1+|\r|)^{1+\nu'}}\p_q \r\ln\lt|\ds\f{1+t}{1+|\r|}\rt|.\eqno{(5.13)}$$
Note that $\p_r\r=\p_p\r+\p_q\r=(1-D_{LL}/4-E_{LL}/4)\p_q\r$ due to $L^i_2\p_i\r=0$, then we have by (5.13)
$$|\p_{\r}\p_q \r|\le \ds\f{C\ve}{(1+|\r|)^{1+\nu'}}(1-D_{LL}/4-E_{LL}/4)^{-1}\p_r \r\ln\lt|\ds\f{1+t}{1+|\r|}\rt|.\eqno{(5.14)}$$
In addition, by a direct computation, we have
\begin{align*}
\p_\r\p_q\r&=\p_\r\lt[(1-\ds\f{1}{4}D_{LL}-\ds\f{1}{4}E_{LL})^{-1}\p_r\r\rt]\\
&=(1-\ds\f{1}{4}D_{LL}-\ds\f{1}{4}E_{LL})^{-1}\p_\r\p_r\r+\ds\f{1}{4}(1-\ds\f{1}{4}D_{LL}-\ds\f{1}{4}E_{LL})^{-2}\p_\r(D_{LL}+E_{LL})\p_r\r.
\end{align*}
Thus it follows that
\begin{align*}
|\p_{\r}\p_r\r|&\le\lt|(1-\ds\f{1}{4}D_{LL}-\ds\f{1}{4}E_{LL})\p_\r\p_q\r\rt|
+\ds\f{1}{4}\lt|(1-\ds\f{1}{4}D_{LL}-\ds\f{1}{4}E_{LL})^{-1}\p_\r(D_{LL}+E_{LL})\p_r\r\rt|\\
&\le\ds\f{C\ve}{(1+|\r|)^{1+\nu'}}(\p_r\r)\ln\lt|\ds\f{1+t}{1+|\r|}\rt|+\ds\f{1}{2}\lt|\r^{-1}_q\p_q(D_{LL}+E_{LL})\rt|\p_r\r\\
&\le\ds\f{C\ve}{(1+|\r|)^{1+\nu'}}(\p_r\r)\ln\lt|\ds\f{1+t}{1+|\r|}\rt|
+\lt(\ds\f{1+t}{1+|\r|}\rt)^{C\ve}\ds\f{C\ve}{(1+t)(1+|\r|)^{\nu'}}\p_r\r.\tag{5.15}
\end{align*}
Without loss of generality, $t>2$ is assumed. When $(t,x)\in H$,  $|\r|\le t/2$,
and we further have $\lt(\ds\f{1+t}{1+|\r|}\rt)^{C\ve}< \ds\f{1+t}{1+|\r|}$ and $\ln\ds\f{1+t}{1+|\r|}
\ge \ln\ds\f{3}{2}$. This, together with (5.15), yields
\begin{align*}
|\p_{\r}\p_r\r|&\le\ds\f{C\ve}{(1+|\r|)^{1+\nu'}}(\p_r\r)\ln\lt|\ds\f{1+t}{1+|\r|}\rt|+\ds\f{C\ve}{(1+|\r|)^{1+\nu'}}\p_r\r\\
&\le \ds\f{C\ve}{(1+|\r|)^{1+\nu'}}\p_r\r\ln\lt|\ds\f{1+t}{1+|\r|}\rt|\\
&\le \ds\f{C\ve\ln|1+t|}{(1+|\r|)^{1+\nu'}}\p_r\r.\tag{5.16}
\end{align*}
When $(t,x)\not\in H$, $\r= r-t$ and  $\p_\r\p_r\r=0$. Therefore, (5.1) of Lemma 5.1 holds.

Based on the preparations above, we now start to show Proposition 5.2. Using $\p^{\text{\bf k}} Z^I u$
$(0\le\text{\bf k} \le k)$ instead of $u$ in (4.5) yields
\begin{align*}
&\ds\int_{\Sigma_t}|\p\p^{\text{\bf k}} Z^I u|^2 W dx +\int^t_0\int_{\Sigma_\tau}(1+|q|)^{-3/2}|\bar Z \p^{\text{\bf k}} Z^I u|^2 W dxd\tau\\
&\le C\int_{\Sigma_0}|\p \p^{\text{\bf k}} Z^I u|^2 W dx
+\int^t_0\int_{\Sigma_\tau}\ds\f{C\ve}{1+\tau}|\p\p^{\text{\bf k}} Z^I u|^2 W dxd\tau\\
&\qquad +\int^t_0\int_{\Sigma_\tau}\ds\f{C(1+\tau)}{\ve}|\wt\sq_g \p^{\text{\bf k}} Z^I u|^2 W dxd\tau.\tag{5.17}
\end{align*}
Thus, we have
\begin{align*}
&E_{k,i}(t)+\ds\sum_{{\text{\bf k}}\le k,I\le i}\ds\int^t_0\int_{\Sigma_\tau}(1+|q|)^{-3/2}|\bar Z \p^{\text{\bf k}} Z^I u|^2 W dxd\tau\\
&\quad \le C E_{k,i}(0)+\int^t_0 \ds\f{C\ve}{1+\tau}E_{k,i}(\tau)d\tau+\ds\sum_{{\text{\bf k}}\le k, I\le i}\int^t_0
\int_{\Sigma_\tau}\ds\f{C(1+\tau)}{\ve}|\wt\sq_g \p^{\text{\bf k}} Z^I u|^2 W dxd\tau.\tag{5.18}
\end{align*}

Starting from (5.18), we now show (5.10). The proof procedure will be divided into
the following three cases.

\vskip 0.2 true cm
{\bf Case 1. \;\;$i=0$}

Without loss of generality, we may assume $1\le k\le N$ since (5.10) obviously holds for $k=0$. It follows from $\wt\sq_g u=0$ and (3.34) that
for $\text{\bf k}\le k$
$$\wt\sq_g \p^{\text{\bf k}}u=\ds\sum_{s+n={\text{\bf k}}, n<{\text{\bf k}}}C^{lm}_{ ij}(\p^s  D^{ij})\p_l\p_m \p^n u+\ds\sum_{s+n={\text{\bf k}}, n<{\text{\bf k}}}G(\p^s u,\p^n u).\eqno{(5.19)}$$
We note that the first term $\ds\sum_{s+n={\text{\bf k}}, n<{\text{\bf k}}}C^{lm}_{ ij}(\p^s  D^{ij})\p_l\p_m \p^n u$ in (5.19)
have been treated in Lemma 9.2 of [20] as follows
$$\lt|\ds\sum_{s+n={\text{\bf k}},n<{\text{\bf k}}}C^{lm}_{ ij}(\p^s  D^{ij})\p_l\p_m \p^n  u\rt| \le \ds\f{C\ve}{1+t}\ds\sum_{n={\text{\bf k}}+1}|\p^n u|+\ds\f{C\ve}{(1+t)^{1-C\ve}}\ds\sum_{1\le n\le {\text{\bf k}}}|\p^n u|.\eqno{(5.20)}$$
Therefore, we just need to treat the second term in (5.19). We rewrite it as
$$\ds\sum_{s+n={\text{\bf k}},n<{\text{\bf k}}}G(\p^s u,\p^n u)=I_1+I_2\eqno{(5.21)}$$
with
\begin{align*}
&I_1=\ds\sum_{n={\text{\bf k}}-1} G(\p u,\p^n u)+G(\p^{\text{\bf k}} u, u),\\
&I_2=\ds\sum_{s+n={\text{\bf k}},2\le s\le {\text{\bf k}}-1,1\le n\le {\text{\bf k}}-2}G(\p^s u,\p^n u).
\end{align*}
By Lemma 2.5, one has
$$
I_1\le C\big(|\bar Z\p u|\ds\sum_{n={\text{\bf k}}+1} |\p^n u|+|\p^2 u||\bar Z\p^{\text{\bf k}} u|\big).\eqno{(5.22)}
$$
For the first term on the right hand side of (5.22),  we use $|\bar Z \p u|\le C(1+t)^{-1}|Z\p u|$. And for the second term we use $|\p^2 u|\le C\ve(1+t)^{-1+C\ve}(1+|q|)^{-1-\nu'}$.
Then by (5.9) and  $\nu'>1/2$  we obtain
\begin{align*}
|I_1|\le &C(1+t)^{-1}|Z\p u|\ds\sum_{n={\text{\bf k}}+1} |\p^n u|+C\ve(1+t)^{-1+C\ve}(1+|q|)^{-1-\nu'}|\bar Z\p^{\text{\bf k}} u|\\
\le &C(1+t)^{-1}\ds\f{C\ve}{1+t}(1+t)^{C\ve}\ds\sum_{n={\text{\bf k}}+1} |\p^n u|
+C\ve(1+t)^{-1+C\ve}(1+|q|)^{-3/2}|\bar Z\p^{\text{\bf k}} u|.\tag{5.23}
\end{align*}
In addition, by (5.9) for  $s< N'$ or $n+1<N'$, we have
$$
|I_2|\le C\ds\sum_{m+l={\text{\bf k}}+3, 2\le m\le {\text{\bf k}}, 1\le l\le {\text{\bf k}}}|\p^m u||\p^l u|
\le \ds\f{C\ve}{(1+t)^{1-C\ve}}\ds\sum_{1\le n\le {\text{\bf k}}}|\p^n u|.\eqno{(5.24)}
$$
Then it follows (5.23)-(5.24) that
\begin{align*}
&\ds\sum_{\text{\bf k}\le k}\int^t_0\int_{\Sigma_\tau}\ds\f{(1+\tau)}{\ve}\bigl|\ds\sum_{s+n={\text{\bf k}}, n<{\text{\bf k}}}G(\p^s u,\p^n u)\bigr|^2 W dxd\tau\\
\le &\int^t_0\ds\f{C\ve}{1+\tau}E_{k,0}(\tau)d\tau+\int^t_0\ds\f{C\ve}{(1+\tau)^{1-C\ve}}E_{k-1,0}(\tau)d\tau\\
&\quad +\ds\sum_{\text{\bf k}\le k}\int^t_0\int_{\Sigma_\tau}C\ve(1+t)^{-1+C\ve}(1+|q|)^{-3/2}|\bar Z \p^{\text{\bf k}} u|^2 W dxd\tau.\tag{5.25}
\end{align*}
Combining (5.25) and (5.18) yields
\begin{align*}
E_{k,0}(t)&<E_{k,0}(t)+\ds\sum_{{\text{\bf k}}\le k}\int^t_0\int_{\Sigma_\tau}[1-C\ve(1+t)^{-1+C\ve}](1+|q|)^{-3/2}|\bar Z \p^{\text{\bf k}} u|^2 Wdxd\tau\\
&\le C E_{k,0}(0)+\int^t_0\ds\f{C\ve}{1+\tau}E_{k,0}(\tau)d\tau+\int^t_0\ds\f{C\ve}{(1+\tau)^{1-C\ve}}E_{k-1,0}(\tau)d\tau,\tag{5.26}
\end{align*}
which completes the proof of (5.10) for $i=0$.
\vskip 0.2 true cm
{\bf Case 2. \;\;$k=0$}

By (3.33) and $\wt\sq_g u=0$, we have
$$\wt\sq_g Z^I u=\ds\sum_{J+K\le I,K<I}C^{\;I \;\;lm}_{JK ij}(Z^J D^{ij})\p_l\p_m Z^K u
+\ds\sum_{J+K\le I,K<I}G(Z^J u,Z^K u).\eqno{(5.27)}$$
For the first term $\ds\sum_{J+K\le I,K<I}C^{\;I \;\;lm}_{JK ij}(Z^J D^{ij})\p_l\p_m Z^K u$  in (5.27) has
been treated  in Lemma 9.2 of [20]) as follows

\begin{align*}
&\ds\f{1}{\ve}\ds\sum_{I\le i}\ds\int^t_0 (1+\tau)\int\lt|\ds\sum_{J+K\le I,K<I}C^{\;I \;\;lm}_{JK ij}(Z^J D^{ij})\p_l\p_m Z^K u\rt|^2 W dx\\
&\le \int^t_0\ds\f{C\ve}{1+\tau}E_{0,i}(\tau)d\tau+\int^t_0\ds\f{C\ve}{(1+\tau)^{1-C\ve}}E_{1,i-1}(\tau)d\tau.\tag{5.28}
\end{align*}
Hence we only deal with the second term $\ds\sum_{J+K\le I,K<I}G(Z^J u,Z^K u)$ in (5.27). Note that
\begin{align*}
&\ds\sum_{J+K\le I,K<I}G(Z^J u,Z^K u)=G(Z^I u,u)+\ds\sum_{J\le 1, K=I-1}
G(Z^J u,Z^K u)\\
&\qquad \qquad \qquad \qquad \qquad +\ds\sum_{J\le I-1,K\le I-2}G(Z^J u, Z^K u).\tag{5.29}
\end{align*}
For term $G(Z^I u,u)$, as treated for $I_1$ in (5.22), we have
\begin{align*}
|G(Z^I u,u)|&\le C\lt(|\bar Z Z^I u||\p^2 u|+|\p Z^I u||\bar Z \p u|\rt)\\
&\le C\ve(1+t)^{-1+C\ve}(1+|q|)^{-3/2}|\bar Z Z^I u|+\ds\f{C\ve}{(1+t)^2}(1+t)^{C\ve}|\p Z^I u|.
\end{align*}
For the left terms $\ds\sum_{J\le 1, K=I-1}G(Z^J u,Z^K u)$  and  $\ds\sum_{J\le I-1,K\le I-2}G(Z^J u,Z^K u)$
in (5.29), their estimates follow directly from (5.9).

Then proceeding as in Case 1, we finally arrive at
$$E_{0,i}(t)\le CE_{0,i}(0)+\int^t_0\ds\f{C\ve}{1+\tau}E_{0,i}(\tau)d\tau
+\int^t_0\ds\f{C\ve}{(1+\tau)^{1-C\ve}}E_{1,i-1}(\tau)d\tau.$$
\vskip 0.2 true cm
{\bf Case 3. \;\;$k\ge1$, $i\ge 1$}

The case proceeds similarly as in Case 1 and 2.

Combining all the three cases above, we complete the proof of Proposition 5.2.\qquad $\square$

Based on Proposition 5.2, as in the proof of Proposition 9.1 of [20], we have

\vskip 0.3 true cm
{\bf Proposition 5.3.} {\it Assume that the assumptions in Proposition 5.2 are valid. Then for $0\le t< T$,
$$E_{k,i}(t)\le C\ds\sum^i_{l=0} E_{k+l,i-l}(0)(1+t)^{C\ve},     \;\;\;\;k+i\le N.  \eqno{(5.30)}$$
where $C>0$ is independent of $T$.}

{\bf Proof.} Since we have established the crucial Proposition 5.2,
the proof of (5.30) is completely similar to that in Proposition 9.1 of [20] by
Gronwall's inequality,
we omit the proof here.\qquad $\square$
\vskip 0.4 true cm
\centerline{\bf $\S 6$. The proof of Theorem 1.1.}
\vskip 0.4 true cm
In this section, we complete the proofs of the weak decay estimate (3.1) and Theorem 1.1 by continuity method.
Let $N\ge 8$ be fixed and set
$$E_N (t)=\ds\sum_{I\le N}\int |\p Z^I u(t,x)|^2 dx.$$
Without loss of generality, we assume
$$E_N (0)\le \ve^2$$
and also assume for fixed $0<\dl<\ds\f12$ and $0\le t\le T$,
$$E_N (t)\le C\ve^2(1+t)^\dl.\eqno{(6.1)}$$

Next, we derive the weak decay estimate (3.1). First it follows from (3.33) and $\wt\sq_g u=0$ that
\begin{align*}
|\sq Z^I u|&=|(\wt\sq_g-D^{ij}\p_i\p_j-E^{ij}\p_i\p_j)Z^I u|\\
&\le |\wt\sq_g Z^I u|+C|u  \;{\p} ^2 Z^I u|+|G(u,Z^I u)|\\
&\le C\ds\sum_{J+K\le I}|Z^J u||\p^2 Z^K u|+ \ds\sum_{J+K\le I}|G(Z^J u,Z^K u)|\\
&\le \ds\f{C}{1+|q|}\ds\sum_{J+K\le I+1}|Z^J u||\p Z^K u|+C\ds\sum_{J+K\le I}|\p Z^J u||{\p}^2 Z^K u|.\tag{6.2}
\end{align*}
Then by H\"older inequality and $\|(1+|q|)^{-1}Z^J u\|_{L^2}\le\|\p Z^J u\|_{L^2}$, we have for $I\le N-3$
\begin{align*}
&\ds\sum_{L\le 2}\int |Z^L\sq Z^I u(t,x)|dx\\
&\le C\ds\sum_{L\le I+2}\int |\sq Z^L u(t,x)|dx\\
&\le C\ds\sum_{J+K\le I+3}\int \ds\f{1}{1+|q|}|Z^J u||\p Z^K u|dx+C\ds\sum_{J+K\le I+2}\int |\p Z^J u||{\p}^2 Z^K u|dx\\
&\le C\ds\sum_{J+K\le I+3}\|\p Z^J u\|_{L^2}\|\p Z^K u\|_{L^2}\\
&\le C E_N(t).\tag{6.3}
\end{align*}
In addition, by Corollary 10.3 of [20], for  $\phi(0,x)=\p_t\phi(0,x)=0$ when $|x|\ge 1$, then
$$|\phi(t,x)|(1+t+|x|)\le C\ds\sum_{I\le 2}\int^t_0\int\ds\f{|(Z^I\sq\phi)(\tau,y)|}{1+\tau+|y|}dyd\tau+C\ds\sum_{I\le 2}\|\p Z^I \phi(0,.)\|_{L^2}.\eqno{(6.4)}$$
Therefore, it follows from (6.3)-(6.4) together with (6.1) that
\begin{align*}
|Z^I u(t,x)|(1+t+r)&\le C\int^t_0 \ds\f{E_N(\tau)}{1+\tau}d\tau +CE_N(0)\\
&\le C\int^t_0\ds\f{\ve^2(1+\tau)^{\dl} }{1+\tau}d\tau +C\ve^2\\
&\le C\ve^2 \dl^{-1} (1+t)^{\dl} \\
&\le C\ve (1+t)^{\dl}.
\end{align*}
Hence, we obtain (3.1) with $\nu=1-\dl>1/2$, then the weak decay (3.1) holds.
Thus, the estimates in
Proposition 3.9 and Lemma 2.2- Lemma 2.4 are also true. It then follows that all the assumptions of Proposition 5.2 hold.
Then  by (5.30) and  $W\ge 1$, we have
$$E_N(t)\le E_{0,N}(t)\le C\ds\sum^N_{l=0}E_{l,N-l}(0)(1+t)^{C\ve}\le C\ve^2(1+t)^{C\ve},$$
here we use the fact of  $E_{k,i}(0)\le C E_{k+i}(0)$ due to $0<\vp(q)=\ds\int^q_{-\infty}\vp'(s)ds=\ds\int^q_{-\infty}(1+|s|)^{-3/2}ds\le C$.
This yields the proof the Theorem 1.1
by the continuity method and the local existence of the solution to (1.3) (In fact,
only by a rough estimate as in [4, Theorem 1], one can then derive that the $C^{\infty}-$solution $u$ of (1.3)
exists for $t\in [0, T]$ with $T\ge e^{\f{C}{\ve}}$).
In addition,  if we choose $\dl=C\ve$ for some suitable positive $C>0$,  as is seen from the proofs above, we can derive (3.1) with  $\nu=1-C\ve$, which means that
the solution $u$ to (1.3) does not behave like a solution to the 3-D free wave equation.

\vskip 0.4 true cm

{\bf Acknowledgements.} {\it Liu Yingbo and  Yin Huicheng wish to express their gratitude  to
Professor Witt Ingo, University of G\"ottingen for many fruitful discussions in this problem
when they were visiting the Mathematical Institute of the University of G\"ottingen from
February to March of 2013. Ding Bingbing is also thankful to Professor Witt Ingo for his much guidance
and huge helps when she read PhD  in University of G\"ottingen from Octobor of 2012 to February of 2014 under the supervision
of Professor  Witt Ingo.}

\end{document}